\documentclass[11pt,amssymb]{amsart}
\usepackage{amsmath}
\usepackage{latexsym}
\usepackage{amssymb}
\usepackage{amscd}
\input epsf

\usepackage{amsfonts}
\usepackage[all]{xy}

\newcounter{SNO}
\newcounter{THNO}[section]
\renewcommand{\theSNO}{\arabic{section}.\arabic{THNO}.\arabic{SNO}}
\renewcommand{\theTHNO}{\arabic{section}.\arabic{THNO}. }
\def\th#1{
\pagebreak[2]
\refstepcounter{THNO}\par
\vspace{1.5ex}
\par\noindent\begingroup
\leftskip=0em\hspace{0em}{\bf  \theTHNO  #1 } \it }
\def\eth{\par\vspace{1.0ex}\endgroup}
\def\pr{\par\noindent{\bf\ Proof. }}
\def\raz{\refstepcounter{THNO}\par\vspace{1.5ex}
\par\noindent{\theTHNO}\quad}
\def\praz{\refstepcounter{SNO}\par\vspace{1ex}
\par\noindent{\theSNO}\quad}
\def\da{\big\downarrow}
\def\bda{\Big\downarrow}

\def\ss#1{\scriptsize{$#1$}}
\def\ot{\otimes}

\newcounter{INO}
\renewcommand{\theINO}{A.\arabic{INO}}
\def\sec#1{\par\vspace{0.9cm}
\par\noindent
{\large\bf \ #1    }\nopagebreak\par\vspace{0.3cm}}

\def\ap#1{\refstepcounter{INO}\par\vspace{0.5cm}
\par\noindent\begingroup \it
\leftskip=0em\hspace{0em}{\bf\ #1 \theINO\ }}
\def\eap{\par\endgroup}

\def\db#1{ D^b ({#1})}
\def\h#1,#2{{\mathrm{H}\mathrm{o}\mathrm{m}}({#1}\:,\; {#2})}
\def\H#1,#2,#3,#4{{\mathrm{H}\mathrm{o}\mathrm{m}}^{#1}_{#2}({#3}\:,\; {#4})}
\def\E#1,#2,#3,#4{{{\mr E}{\mr x}{\mr t}}^{#1}_{#2}({#3}\:,\; {#4})}
\def\lto{\longrightarrow}

\def\bt{\boxtimes}
\def\cial{\circlearrowleft}
\def\dda{\begin{picture}(4,15)\multiput(2,9)(0,-2){5}{\circle*{1}}
\put(2,1){\vector(0,-1){4}}\end{picture}}
\def\mr#1{\mathrm{#1}}
\def\bdot{\begin{picture}(4,4)\put(2,3){\circle*{1.5}}\end{picture}}
\title{Equivalences of derived categories and K3 surfaces}
\author{Dmitri Orlov}
\date{}
\address{Algebra section,
Steklov Math. Institute, Vavilova 42, 117333 Moscow, RUSSIA  }
\email{orlov@mi.ras.ru}
\begin{document}
\maketitle

\begin{abstract}
We consider derived categories of coherent sheaves on smooth projective
varieties. We prove that any equivalence between them can be represented
by an object on the product. Using this, we give a necessary and sufficient
condition for equivalence of derived categories of two K3 surfaces.
\end{abstract}
%\newpage

%\input algfonts
%\maketitle

%%%%%%%%%%%%%%%%%%%%%%%%%%%%%%%%%%%%%%%%%%%%%%%%%%%%%%%%%%%%%%%%%%
\section*{Introduction}

Let $\db{X}$ be the bounded derived category of coherent sheaves on a smooth
projective variety $X.$ The category $\db{X}$ has the structure of a triangulated
category (see \cite{Ver}, \cite{GM}). We shall consider $\db{X}$
as a triangulated category.

In this paper we are concerned with the problem of description for varieties,
which
have equivalent derived categories of coherent sheaves.

In the paper \cite{Mu1}, Mukai showed that for an abelian variety $A$ and its dual
$\hat{A}$ the derived categories $\db{A}$
and $\db{\hat{A}}$ are equivalent . Equivalences of another type appeared in \cite{BO}.
They
are induced by certain birational transformations which are called flops.

Further, it was proved  in the  paper \cite{BOr} that if
$X$ is a smooth projective variety with  either ample canonical or ample
anticanonical
sheaf, then any other algebraic variety $X'$ such that $\db{X'}\simeq\db{X}$
is biregularly isomorphic to $X.$

The aim of this paper is to give some description for  equivalences  between
derived categories of coherent sheaves. The main result is Theorem \ref{main} of
 $\S 2.$
It says that any full and faithful exact functor $F: \db{M}\lto \db{X}$ having left (or right) adjoint functor can
be represented by
an object $E\in\db{M\times X},$ i.e. $F(\bdot)\cong
R^{\bdot}\pi_*(E\stackrel{L}{\ot}p^*(\bdot)),$ where $\pi$ and $p$ are the
projections on $M$ and $X$ respectively.

In $\S 3,$ basing on the Mukai's results \cite{Mu}, we show that two K3
surfaces $S_1$ and $S_2$ over field $\Bbb{C}$ have equivalent derived categories
of coherent sheaves iff the lattices of transcendental cycles $T_{S_1}$ and
$T_{S_2}$ are Hodge isometric.

I would like to thank A.~Polishchuk for useful notices.

\section{Preliminaries}

\raz We collect here some facts relating to triangulated categories. Recall
that a triangulated category is an additive category with additional structures:

a) {\it an  additive autoequivalence $T : {\mathcal D}\lto {\mathcal
D},$ which is called a translation functor} (we usually write $X[n]$
instead of $T^n(X)$ and $f[n]$ instead of $T^n(f)$),

b) {\it a class of distinguished triangles:}
$$
X\stackrel{u}{\to}Y\stackrel{v}{\to}Z\stackrel{w}{\to}X[1].
$$
And these structures must satisfy the usual set of axioms (see \cite{Ver}).

If $X,$ $ Y$ are objects of a triangulated category ${\mathcal D},$
then ${\H i,{\mathcal D}, X, Y}$ means ${\H{},{\mathcal D}, X,
{Y[i]}}.$

An additive functor $F : {\mathcal D}\lto{\mathcal D}'$ between two
triangulated categories ${\mathcal D}$ and ${\mathcal D}'$ is called
{\sf exact} if

a)  {\it  it commutes with the translation functor, i.e there is fixed
an isomorphism of functors:}
$$
t_F : F\circ T\stackrel{\sim}{\lto}T'\circ F,
$$

b) {\it it takes every distinguished triangle to a distinguished triangle}
(using the isomorphism $t_F,$ we replace $F(X[1])$ by $F(X)[1]$).

The following lemma will be needed for the sequel.
\th{Lemma}\cite{BK} If a functor $G : {\mathcal D}'\lto{\mathcal D}$
is a left (or right) adjoint to an exact functor $F : {\mathcal
D}\lto{\mathcal D}'$ then functor $G$ is also exact . \eth \pr Since
$G$ is the left adjoint functor to $F,$ there exist canonical
morphisms of functors $id_{\mathcal D'}\to F\circ G,\; G\circ F\lto
id_{\mathcal D}.$ Let us consider  the following sequence of natural
morphisms:
$$
G\circ T'\lto G\circ T'\circ F\circ G\stackrel{\sim}{\lto}
G\circ F\circ T\circ G\lto T\circ G
$$
We obtain the natural morphism $G\circ T'\lto T\circ G.$ This
morphism is an isomorphism. Indeed, for any two objects $A\in
{\mathcal D}$ and $B\in {\mathcal D}'$ we have isomorphisms :
$$
\begin{array}{l}
{\h G(B[1]), A}\cong{\h B[1], {F(A)}}\cong{\h B, {F(A)[-1]}}\cong\\\\
{\h B, {F(A[-1])}}\cong{\h G(B), {A[-1]}}\cong{\h G(B)[1], A}\\
\end{array}
$$
This implies that the natural morphism $G\circ T'\lto T\circ G$ is an isomorphism.

Let now $A\stackrel{\alpha}{\lto}B\lto C\lto A[1]$ be a
distinguished triangle in ${\mathcal D}'.$ We have to show that $G$
takes this triangle to a distinguished one.

Let us include the morphism $G(\alpha) : G(A)\to G(B)$ into a distinguished
triangle:
$$
G(A)\lto G(B)\lto Z\lto G(A)[1].
$$
Applying functor $F$ to it, we obtain a distinguished triangle:
$$
FG(A)\lto FG(B)\lto F(Z)\lto FG(A)[1]
$$
(we use the commutation isomorphisms like $T'\circ F\stackrel{\sim}{\to}
F\circ T$ with no mention).

Using morphism $id\to F\circ G,$ we get a commutative diagram:
$$
\begin{array}{ccccccc}
A&\stackrel{\alpha}{\lto}&B&\lto& C&\lto& A[1]\\
\da&&\da&&&&\da\\
FG(A)&\stackrel{FG(\alpha)}{\lto}& FG(B)&\lto& F(Z)&\lto& FG(A)[1]
\end{array}
$$
By axioms of triangulated categories there exists a morphism $\mu : C\to
F(Z)$ that completes this commutative diagram. Since $G$ is left adjoint to
$F,$ the morphism $\mu$ defines $\nu : G(C)\to Z.$ It is clear that
$\nu$ makes  the following diagram commutative:
$$
\begin{array}{ccccccc}
G(A)&{\lto}&G(B)&\lto&G(C)&\lto&G(A)[1]\\
\wr\da&&\wr\da&&\da\rlap{\ss{\nu}}&&\wr\bda\\
G(A)&\lto&G(B)&\lto& Z&\lto&G(A)[1]\\
\end{array}
$$

To prove the lemma, it suffices to show that $\nu$ is an isomorphism.
 For any object $Y\in{\mathcal D}$ let us consider the diagram for
$\mbox{Hom}$:
$$
\begin{array}{cccccccc}
\to{\h G(A)[1], Y}&\to&{\h Z, Y}&\to&{\h G(B), Y}\to\\
\da\wr&&\da\rlap{\ss{\mr{H}_Y (\nu)}}&&
\da\wr\\
\to{\h G(A)[1], Y}&{\to}&{\h G(C), Y}&\to&{\h G(B), Y}
\to\\
\da\wr&&\da\wr&&\da\wr\\
\to{\h A[1], {F(Y)}}&\to&{\h C, {F(Y)}}&\to&{\h B, {F(Y)}}
\to\\
\end{array}
$$
Since the lower sequence is exact, the middle sequence is exact also.
By the lemma about five homomorphisms, for any $Y$ the morphism $\mr{H}(\nu)$ is an
isomorphism . Thus
$\nu : G(C)\to Z$ is an isomorphism too. This concludes the proof. $\Box$

\raz
Let $X^{\bdot} = \{ X^{c}\stackrel{d^{c}}{\to}X^{c+1}\stackrel{d^{c+1}}{\to}
\cdots\to X^0\}$
 be a bounded complex over a triangulated category ${\mathcal D},$
i.e.  all compositions $d^{i+1}\circ d^i$ are equal to $0$ ($c< 0$).

A left Postnikov system, attached to $X^{\bdot},$ is, by definition, a diagram

\begin{picture}(400,100)
\put(24,75){\vector(1,-2){30}}
\put(64,15){\vector(1,2){30}}
\put(104,75){\vector(1,-2){30}}
\put(144,15){\vector(1,2){30}}
\put(184,75){\vector(1,-2){30}}
\put(304,15){\vector(1,2){30}}
\put(344,75){\vector(1,-2){30}}
\put(40,82){\vector(1,0){40}}
\put(120,82){\vector(1,0){40}}
\put(120,2){\vector(-1,0){40}}
\put(200,2){\vector(-1,0){40}}
\put(250,2){\vector(-1,0){10}}
\put(285,2){\vector(-1,0){10}}
\put(360,2){\vector(-1,0){40}}
\put(10,80){$X^{c}$}
\put(90,80){$X^{c+1}$}
\put(170,80){$X^{c+2}$}
\put(32,0){$Y^{c}=X^{c}$}
\put(130,0){$Y^{c+1}$}
\put(210,0){$Y^{c+2}$}
\put(330,80){$X^0$}
\put(370,0){$Y^0$}
\put(290,0){$Y^{-1}$}
\put(252,0){$\cdots$}
\put(3,32){$i_c=id$}
\put(62,43){$ j_c$}
\put(100,32){$ i_{c+1}$}
\put(134,43){$ j_{c+1}$}
\put(178,32){$ i_{c+2}$}
\put(302,43){$ j_{-1}$}
\put(342,32){$ i_0$}
\put(55,85){$ d^{c}$}
\put(135,85){$ d^{c+1}$}
\put(50,60){$ \cial$}
\put(90,26){$ \star$}
\put(130,60){$ \cial$}
\put(170,26){$ \star$}
\put(330,26){$ \star$}
\put(90,7){$ [1]$}
\put(170,7){$ [1]$}
\put(330,7){$ [1]$}
\end{picture}

\bigskip
\noindent in which all triangles marked with $\star$ are
distinguished and triangles marked with $\cial$ are commutative
(i.e. $j_k\circ i_k = d^{k}$). An object $E\in\mbox{Ob}{\mathcal D}$
is called a left convolution of $X^{\bdot},$ if there exists a left
Postnikov system, attached to $X^{\bdot}$ such that $E=Y^0.$ The
class of all convolutions of $X^{\bdot}$ will be denoted by
$\mbox{Tot} (X^{\bdot}).$ Clearly the Postnikov systems and
convolutions are stable under exact functors between triangulated
categories.

The class $\mbox{Tot}(X^{\bdot})$ may contain many non-isomorphic elements
 and may be empty. Further we shall give a sufficient condition for
$\mbox{Tot}(X^{\bdot})$ to be non-empty and for its objects to be isomorphic.
The following lemma is needed for the sequel(see \cite{BBD}).
\th{Lemma}\label{tr}
Let $g$ be a morphism between two objects $Y$ and $Y',$ which  are included into
two distinguished triangles:
$$
\begin{array}{ccccccc}
X&\stackrel{u}{\lto}&Y&\stackrel{v}{\lto}&Z&\stackrel{w}{\lto}&X[1]\\
\dda\rlap{\ss{f}}&&\da\rlap{\ss{g}}&&\dda\rlap{\ss{h}}&&\dda\rlap{\ss{f[1]}}\\
X'&\stackrel{u'}{\lto}&Y'&\stackrel{v'}{\lto}&Z'&\stackrel{w'}{\lto}&X'[1]
\end{array}
$$
If $v'gu=0,$ then there exist morphisms $f : X\to X'$ and $h : Z\to Z'$ such
that the triple $(f, g, h)$ is a morphism of triangles.

If, in addition, ${\h X[1], {Z'}}=0$ then this triple is uniquely determined
by $g.$
\eth
Now we prove two lemmas which are generalizations of the previous one for
Postnikov diagrams.
\th{Lemma}\label{pd1}
 Let $X^{\bdot} = \{ X^{c}\stackrel{d^{c}}{\to}X^{c+1}\stackrel{d^{c+1}}{\to}
\cdots\to X^0\}$
 be a bounded complex over a triangulated category ${\mathcal D}.$
Suppose it satisfies the following condition:
\begin{equation}\label{ex}
{\H i, {}, X^a, {X^b}}=0\; \mbox{ for }\; i<0 \;\mbox{ and }\; a<b.
\end{equation}
  Then there exists a convolution for this complex and all convolutions are
isomorphic (noncanonically).

 If, in addition,
\begin{equation}\label{une}
{\H i, {}, X^a, {Y^0}}=0 \;\mbox{ for }\; i<0\; \mbox{ and for all }\; a
\end{equation}
 for some convolution
$Y^0$ (and, consequently, for any one), then all convolutions are canonically
isomorphic.
\eth
\th{Lemma}\label{pd2}
 Let $X^{\bdot}_1$ and $X^{\bdot}_2$ be  bounded complexes that
satisfy (\ref{ex}),
and let $(f_c,...,f_0)$ be
a morphism of these complexes:
$$
\begin{array}{ccccccc}
X^c_1&\stackrel{d^c_1}{\lto}&X^{c+1}_1&\lto&\cdots&\lto&X^0_1\\
\da\rlap{\ss{f_c}}&&\da\rlap{\ss{f_{c+1}}}&&&&\da\rlap{\ss{f_0}}\\
X^c_2&\stackrel{d^c_2}{\lto}&X^{c+1}_2&\lto&\cdots&\lto&X^0_2\\
\end{array}
$$
Suppose that
\begin{equation}\label{exm}
{\H i, {}, X^a_1, {X^b_2}}=0 \;\mbox{ for }\; i<0 \;\mbox{ and }\; a<b.
\end{equation}
 Then  for any convolution $Y^0_1$ of $X^{\bdot}_1$
and for any convolution $Y^0_2$ of $X^{\bdot}_2$ there
exists a morphism $f : Y^0_1\to Y^0_2$ that commutes with  the morphism  $f_0.$
 If, in addition,
\begin{equation} \label{unm}
{\H i, {}, X^a_1, {Y^0_2}}=0 \;\mbox{ for }\; i<0\; \mbox{ and for all }\; a
\end{equation}
then this morphism
is unique.
\eth
\pr
We shall prove  both lemmas together.
 Let $Y^{c+1}$ be a cone of the morphism $d^c$:
$$
X^c\stackrel{d^c}{\lto}X^{c+1}\stackrel{\alpha}{\lto} Y^{c+1}\lto X^c [1]
$$
By assumption $d^{c+1}\circ d^c =0$ and ${\h X^c [1], {X^{c+2}}}=0,$ hence
there exists a unique morphism $\bar{d}^{c+1} : Y^{c+1}\to
X^{c+2}$ such that
$\bar{d}^{c+1}\circ
 \alpha = d^{c+1}.$

Let us consider a composition
$d^{c+2}\circ \bar{d}^{c+1} : Y^{c+1}\to X^{c+3}.$
We know that $d^{c+2}\circ \bar{d}^{c+1}\circ\alpha = d^{c+2}\circ d^{c+1}=0 ,$
and at the same time we have ${\h X^c [1], {X^{c+3}}}=0.$ This implies that
the composition $d^{c+2}\circ \bar{d}^{c+1}$ is equal to $0.$

Moreover, consider the distinguished triangle for $Y^{c+1}.$ It can easily be
 checked that
 ${\H i, {}, Y^{c+1},{ X^b}}=0$ for $i<0$ and $b>c+1.$
Hence the complex
$Y^{c+1}\lto X^{c+2}\lto\cdots\lto X^0$  satisfies the condition (\ref{ex}).
By induction, we can suppose  that it has a convolution. This implies that
 the complex $X^{\bdot}$ has a convolution too.
Thus, the class $\mbox{Tot}(X^{\bdot})$ is non-empty.

Now we shall show that under the conditions (\ref{exm}) any morphism  of complexes
can be extended to a morphism of Postnikov systems.

Let us consider  cones $Y^{c+1}_1$ and $Y^{c+1}_2$ of the morphisms $d^c_1$ and
$d^c_2.$ There exists a morphism $g_{c+1} : Y^{c+1}_1\to Y^{c+1}_2$ such that
one has
the morphism of distinguished triangles:
$$
\begin{array}{ccccccc}
X^c_1&\stackrel{d^c_1}{\lto}&X^{c+1}_1&\stackrel{\alpha}{\lto}&
Y^{c+1}_1&\lto&X^c_1 [1]\\
\da\rlap{\ss{f_c}}&&\da\rlap{\ss{f_{c+1}}}&&\da\rlap{\ss{g_{c+1}}}&
&\da\rlap{\ss{f_c [1]}}\\
X^c_2&\stackrel{d^c_2}{\lto}&X^{c+1}_2&\stackrel{\beta}{\lto}&
Y^{c+1}_2&\lto&X^c_2 [1]\\
\end{array}
$$
As above, there exist uniquely determined morphisms
$ \bar{d}^{c+1}_i : Y^{c+1}_i\to X^{c+2}_i$ for $i=1,2.$
Consider the following diagram:
$$
\begin{array}{ccc}
Y^{c+1}_1&\stackrel{\bar{d}^{c+1}_1}{\lto}&X^{c+2}_1\\
\da\rlap{\ss{g_{c+1}}}&&\da\rlap{\ss{f_{c+2}}}\\
Y^{c+1}_2&\stackrel{\bar{d}^{c+1}_2}{\lto}&X^{c+2}_2
\end{array}
$$
Let us show that this square is commutative. Denote by $h$  the difference
$f_{c+2}\circ \bar{d}^{c+1}_1 - \bar{d}^{c+1}_2\circ g_{c+1}.$ We have
$h\circ \alpha = f_{c+2}\circ {d}^{c+1}_1 - {d}^{c+1}_2\circ f_{c+1} = 0$
and, by assumption, ${\h X^c_1 [1], {X^{c+2}_2}}=0.$ It follows that $h=0.$
 Therefore, we obtain the  morphism of new complexes:
$$
\begin{array}{ccccccc}
Y^{c+1}_1&\stackrel{\bar{d}^{c+1}_1}{\lto}&X^{c+2}_1&\lto&\cdots&\lto&X^0_1\\
\da\rlap{\ss{g_{c+1}}}&&\da\rlap{\ss{f_{c+2}}}&&&&\da\rlap{\ss{f_0}}\\
Y^{c+1}_2&\stackrel{\bar{d}^{c+1}_2}{\lto}&X^{c+2}_2&\lto&\cdots&\lto&X^0_2\\
\end{array}
$$
It can easily be checked that these complexes satisfy the conditions (\ref{ex})
 and (\ref{exm}) of the lemmas.
By the induction hypothesis, this morphism can be extended to a morphism of Postnikov
systems, attached to these complexes. Hence there exists a morphism of
Postnikov systems, attached to $X^{\bdot}_1$ and $X^{\bdot}_2.$

Moreover, we  see that if all morphisms $f_i$ are isomorphisms, then a
morphism of Postnikov systems is an isomorphism too. Therefore, under the
condition (\ref{ex}) all objects from the class $\mbox{Tot}(X^{\bdot})$ are isomorphic.

Now let us  consider a morphism of the rightmost distinguished triangles of
Postnikov systems:
$$
\begin{array}{ccccccc}
Y^{-1}_1&\stackrel{j_{1, -1}}{\lto}&X^0_1&\stackrel{i_{1, 0}}{\lto}&
Y^0_1&\lto&Y^{-1}_1 [1]\\
\da\rlap{\ss{g_{-1}}}&&\da\rlap{\ss{f_0}}&&\da\rlap{\ss{g_0}}&
&\da\rlap{\ss{g_{-1} [1]}}\\
Y^{-1}_2&\stackrel{j_{2, -1}}{\lto}&X^0_2&\stackrel{i_{2, 0}}{\lto}&
Y^0_2&\lto&Y^{-1}_2 [1]\\
\end{array}
$$

If the complexes  $X^{\bdot}_i$ satisfy the condition (\ref{unm}) (
i.e.${\H i, {}, X^a_1, {Y^0_2}}=0$ for $ i<0 $ and all $ a$),
then we get ${\h Y^{-1}_1 [1], {Y^0_2}}=0.$
It follows from Lemma \ref{tr}  that $g_0$ is uniquely determined.
This concludes the proof of  both lemmas. $\Box$

\section{Equivalences of derived categories}

\raz Let $X$ and $M$ be  smooth projective varieties over field $k.$
Denote by $\db{X}$ and $\db{M}$ the bounded derived categories of
coherent sheaves on $X$ and $M$ respectively. Recall that a derived
category has the structure of a triangulated category.

For every object $E\in \db{M\times X}$ we can define an exact functor
$\Phi_E$ from $\db{M}$ to $\db{X}.$ Denote by $p$ and $\pi$ the projections
of ${M\times X}$ onto $M$ and $X$ respectively:
$$
\begin{array}{ccc}
M\times X&\stackrel{\pi}{\lto}&X\\
\llap{\ss{p}}\da&&\\
M&&
\end{array}
$$
Then $\Phi_E$ is defined by the following formula:
\begin{equation}\label{dfun}
\Phi_E(\bdot):= \pi_*(E\otimes p^*(\bdot))
\end{equation}
(we always shall write shortly $f_* , f^*, \otimes$ and etc. instead of $R^{\bdot}f_*,
L^{\bdot}f^*, \stackrel{L}{\otimes},$ because we consider only derived functors).

The functor $\Phi_E$ has the left and the right adjoint functors $\Phi_E^*$ and
$\Phi_E^!$ respectively, defined by the following formulas:
$$
\begin{array}{l}
\Phi_E^* (\bdot) = p_*(E^{\vee}\ot \pi^*(\omega_X [dimX]\ot (\bdot))),\\\\
\Phi_E^! (\bdot) =\omega_M [dimM]\ot  p_*(E^{\vee}\ot (\bdot)),
\end{array}
$$
where $\omega_X$ and $\omega_M$ are the canonical sheaves on $X$ and
$M,$ and $E^{\vee}:= {\pmb R}^{\bdot}{\mathcal H}om( E, {\mathcal
O}_{M\times X} ).$

Let $F$ be an exact functor from the derived category $\db{M}$ to the
derived category $\db{X}.$ Denote by $F^*$ and $F^!$ the left and the right
adjoint functors for $F$ respectively, when they exist. Note that if there
exists the left adjoint functor $F^*,$ then the right adjoint functor $F^!$
 also exists and
$$
F^! = S_M\circ F^*\circ S_X^{-1},
$$
where $S_X$ and $S_M$ are Serre functors on $\db{X}$ and $\db{M}.$ They are
equal to $(\bdot)\ot\omega_X [dimX]$ and $(\bdot)\ot\omega_M [dimM]$ (see
\cite{BK}).

What can we say about the category of all exact functors between $\db{M}$
and $\db{X}$? It seems to be true that any functor can be represented by an
object on the product $M\times X$ for smooth projective varieties $M$ and $X.$
  But we are unable prove it. However, when $F$ is full and faithfull, it can be
represented. The  main result of this chapter is the following theorem.

\th{Theorem}\label{main}
Let $F$ be an exact functor from $\db{M}$ to $\db{X},$ where $M$ and $X$ are
smooth projective varieties. Suppose $F$ is full and faithful and has
the right (and,consequently, the left) adjoint functor.

Then there exists an object $E\in\db{M\times X}$ such
 that $F$ is isomorphic to the functor $\Phi_E$ defined by the rule
(\ref{dfun}), and this object is unique up to isomorphism.
\eth

\raz Let $F$ be an exact functor from a derived category
$\db{\mathcal A}$ to a derived category $\db{\mathcal B}.$
 We say that $F$ is {\sf bounded} if there exist
$z\in {\Bbb{Z}}, n\in \Bbb{N}$ such that for any $A\in{\mathcal A}$
the cohomology objects ${H}^i (F(A))$ are equal to $0$ for $i\not\in
[z, z+n].$

\th{Lemma} Let $M$ and $X$ be smooth projective varieties. If an
exact functor $F : \db{M}\lto\db{X}$ has a left  adjoint functor
then it is bounded. \eth \pr Let $G : \db{X}\lto\db{M}$ be a left
adjoint functor to $F.$ Take a very ample invertible sheaf
${\mathcal L}$ on $X.$ It gives the embedding $i : X \hookrightarrow
{\Bbb{P}}^N.$ For any $i<0$ we have right resolution of the sheaf
${\mathcal O}(i)$ on ${\Bbb{P}}^N$ in terms of the sheaves
${\mathcal O}(j),$
 where $j=0, 1,.., N$(see \cite{Be}). It is easily seen that this resolution
is of the form
$$
{\mathcal O}(i)\stackrel{\sim}{\lto}\Bigl\{ V_0\otimes{\mathcal
O}\lto V_1\otimes{\mathcal O}(1)\lto\cdots\lto V_N\otimes{\mathcal
O}(N)\lto 0 \Bigl\}
$$
where all $V_k$ are vector spaces. The restriction of this
resolution to $X$ gives us the resolution of the sheaf ${\mathcal
L}^{ i}$ in terms of the sheaves ${\mathcal L}^{ j},$ where $j=0,
1,..., N.$ Since the functor $G$ is exact that there exist $z'$ and
$n'$ such that ${H}^{k}(G({\mathcal L}^{ i}))$ are equal $0$ for
$k\not\in [z', z'+n'].$ This follows from the existence of the
spectral sequence
$$
E^{p,q}_1 = V_p\otimes {H}^q(G({\mathcal L}^{p})) \Rightarrow
{H}^{p+q}(G({\mathcal L}^{i})).
$$
As all nonzero terms  of this spectral sequence are concentrated in
some rectangle, so it follows that  for all $i$ cohomologies
${H}^{\bdot}(G({\mathcal L}^{i}))$  are concentrated  in some
segment.

Now, notice that if ${\H j, {}, {\mathcal L}^{i}, {F(A)}}= 0$ for
all $i\ll0,$ then ${H}^j (F(A))$ is equal to $0.$ Further,
 by assumption,
 the functor $G$ is left adjoint to $F,$ hence
$$
{\H j, {}, {\mathcal L}^{i}, {F(A)}}\cong
{\H j, {}, G({\mathcal L}^{i}), A}.\\
$$

If now $A$ is a sheaf on $M,$ then ${\H j, {}, G({\mathcal L}^{i}),
A}=0$ for all $i<0$ and $j\not\in [-z'-n', -z'+ dimM],$ and thus
${H}^j (F(A))=0$ for the same $j.$ $\Box$

\th{Remark}\label{boun} We shall henceforth assume that for any
sheaf ${\mathcal F}$ on $M$ the cohomology objects ${H}^i
(F({\mathcal F}))$ are nonzero only if $i\in [-a, 0].$ \eth

\raz \label{cons} Now we begin  constructing an object $E\in\db{M\times X}.$
Firstly, we shall consider a closed embedding $j : M\hookrightarrow{\Bbb{P}}^N$
and shall construct an object $E'\in \db{{\Bbb{P}}^N \times X}.$ Secondly,
we shall show that there exists an object  $E\in\db{M\times X}$ such that
$E'=(j\times id)_* E.$ After that we shall prove that functors $F$ and $\Phi_E$ are isomorphic.

Let ${\mathcal L}$ be a very ample invertible sheaf on $M$ such that
${\mr H}^i ({\mathcal L}^{k})=0$ for any $k>0,$ when $i\not=0.$ By
$j$ denote the closed embedding $j : M\hookrightarrow{\Bbb{P}}^N$
with respect to ${\mathcal L}.$

Recall that there exists a resolution of the diagonal  on the product
${\Bbb{P}}^N \times {\Bbb{P}}^N$ (see\cite{Be}). Let us consider the following
complex of sheaves on the product:
\begin{equation}\label{di}
0\to {\mathcal O}(-N)\bt\Omega^N (N)\stackrel{d_{-N}}{\to}{\mathcal
O}(-N+1)\bt \Omega^{N-1} (N-1)\to\cdots\to{\mathcal
O}(-1)\bt\Omega^{1} (1) \stackrel{d_{-1}}{\to}{\mathcal
O}\bt{\mathcal O}
\end{equation}
This complex is a resolution of the structure sheaf ${\mathcal
O}_{\Delta}$ of the diagonal $\Delta.$

Now by $F'$ denote the functor from $\db{{\Bbb{P}}^N}$ to $\db{X},$ which
is the composition $F\circ j^*.$ Consider the product
$$
\begin{array}{ccc}
{\Bbb{P}}^N \times X&\stackrel{\pi^{'}}{\lto}&X\\
\llap{\ss{q}}\da&&\\
{\Bbb{P}}^N
\end{array}
$$
Denote by
$$
d'_{-i}\in{\H {}, {{\Bbb{P}}^N \times X}, {\mathcal O}(-i)\bt
F'(\Omega^i (i)), {{\mathcal O}(-i+1)\bt F'(\Omega^{i-1} (i-1))}}
$$
the image $d_{-i}$
under the following through map.
$$
\begin{array}{l}
{\h {\mathcal O}(-i)\bt \Omega^i (i),
{{\mathcal O}(-i+1)\bt \Omega^{i-1} (i-1)}}\stackrel{\sim}{\lto}\\\\
{\h {\mathcal O}\bt \Omega^i (i),
{{\mathcal O}(1)\bt \Omega^{i-1} (i-1)}}\stackrel{\sim}{\lto}\\\\
{\h  \Omega^i (i),
{{\mr H}^0 ({\mathcal O}(1))\ot \Omega^{i-1} (i-1)}}\lto\\\\
{\h  F'(\Omega^i (i)),
{{\mr H}^0 ({\mathcal O}(1))\ot F'(\Omega^{i-1} (i-1))}}\stackrel{\sim}{\lto}\\\\
{\h {\mathcal O}\bt F'(\Omega^i (i)),
{{\mathcal O}(1)\bt F'(\Omega^{i-1} (i-1))}}\stackrel{\sim}{\lto}\\\\
{\h {\mathcal O}(-i)\bt F'(\Omega^i (i)), {{\mathcal O}(-i+1)\bt
F'(\Omega^{i-1} (i-1))}}
\end{array}
$$
It can easily be checked that the composition $d_{-i+1}\circ d_{-i}$ is
equal to $0.$ We omit the check.

Consider the following complex $C^{\bdot}$
$$
C^{\bdot} :=\{ {\mathcal O}(-N)\bt F'(\Omega^N
(N))\stackrel{d'_{-N}}{\lto} \cdots\lto {\mathcal O}(-1)\bt
F'(\Omega^{1} (1))\stackrel{d'_{-1}}{\lto} {\mathcal O}\bt
F'({\mathcal O}) \}
$$
over the derived category $\db{{\Bbb{P}}^N \times X}.$ For $l<0$ we have
$$
\begin{array}{l}
{\H l, {}, {\mathcal O}(-i)\bt F'(\Omega^i (i)),
{{\mathcal O}(-k)\bt F'(\Omega^{k} (k))}}\cong\\\\
{\H l, {}, {\mathcal O}\bt F'(\Omega^i (i)),
{{\mr H}^0 ({\mathcal O}(i-k))\ot F'(\Omega^{k} (k))}}\cong\\\\
{\H l, {}, j^* (\Omega^i (i)), {{\mr H}^0 ({\mathcal O}(i-k))\ot j^*
(\Omega^{k} (k))}}=0
\end{array}
$$
Hence, by Lemma \ref{pd1}, there exists a convolution of the complex
$C^{\bdot},$ and all convolutions are isomorphic. By $E'$ denote
some convolution of $C^{\bdot}$ and by $\gamma_0$ denote the
morphism ${\mathcal O}\bt F'({\mathcal
O})\stackrel{\gamma_0}{\lto}E'.$ (Further we shall see that all
convolutions of $C^{\bdot}$ are canonically isomorphic). Now let
$\Phi_{E'}$ be the functor from $\db{{\Bbb{P}}^N}$ to $\db{X},$
defined by (\ref{dfun}).

\th{Lemma}\label{tran} There exist canonically defined isomorphisms
$f_k : F'({\mathcal O}(k)) \stackrel{\sim}{\lto}\Phi_{E'}({\mathcal
O}(k))$ for all $k\in{\Bbb{Z}},$ and these isomorphisms are
functorial, i.e. for any $\alpha : {\mathcal O}(k)\to {\mathcal
O}(l)$ the following diagram commutes
$$
\begin{array}{ccc}
  F'({\mathcal O}(k))&\stackrel{F'(\alpha)}{\lto}&F'({\mathcal O}(l))\\
\llap{\ss{f_k}}\da&&\da\rlap{\ss{f_l}}\\
\Phi_{E'}({\mathcal
O}(k))&\stackrel{\Phi_{E'}(\alpha)}{\lto}&\Phi_{E'}({\mathcal O}(l))
\end{array}
$$
\eth
\pr
At first, assume that $k\ge0.$

Consider the resolution (\ref{di}) of the diagonal
$\Delta\subset{\Bbb{P}}^N \times{\Bbb{P}}^N$ and, after tensoring it
with ${\mathcal O}(k)\bt{\mathcal O},$ push forward onto the second
component. We get the following resolution of ${\mathcal O}(k)$ on
${\Bbb{P}}^N$
$$
\{ {\mr H}^0 ({\mathcal O}(k-N))\ot\Omega^N (N){\lto}\cdots\lto{\mr
H}^0 ({\mathcal O}(k-1)) \ot\Omega^{1}(1){\lto}{\mr H}^0 ({\mathcal
O}(k))\ot{\mathcal O} \}\stackrel{\delta_k}{\lto}{\mathcal O}(k)
$$
Consequently $F'({\mathcal O}(k))$ is a convolution of the complex
$D^{\bdot}_k$:
$$
{\mr H}^0 ({\mathcal O}(k-N))\ot F'(\Omega^N
(N)){\lto}\cdots\lto{\mr H}^0 ({\mathcal O}(k-1)) \ot
F'(\Omega^{1}(1)){\lto}{\mr H}^0 ({\mathcal O}(k))\ot F'({\mathcal
O})
$$
over $\db{X}.$

On the other hand, let us consider the complex $C^{\bdot}_k := q^*
{\mathcal O}(k)\ot C^{\bdot}$ on ${\Bbb{P}}^N\times X$ with the
morphism $\gamma_k : {\mathcal O}(k)\bt F'({\mathcal O})\lto
q^*{\mathcal O}(k)\ot E',$ and push it forward onto the second
component. It follows from the construction of the complex
$C^{\bdot}$ that $\pi'_* (C^{\bdot}_k)=D^{\bdot}_k.$ So we  see that
$F'({\mathcal O}(k))$ and $\Phi_{E'}({\mathcal O}(k))$ both are
 convolutions of the same complex $D^{\bdot}_k.$

By assumption the functor $F$ is full and faithful, hence, if
${\mathcal G}$ and ${\mathcal H}$ are locally free sheaves on
${\Bbb{P}}^N$  then we have
$$
{\H i, {}, F'({\mathcal G}), {F'({\mathcal H})}}= {\H i, {}, j^*
({\mathcal G}), {j^*({\mathcal H})}}=0
$$
for $i<0.$ Therefore the complex $D^{\bdot}_k$ satisfies the
conditions (\ref{ex}) and (\ref{une}) of Lemma \ref{pd1} Hence
there exists a uniquely defined isomorphism $f_k : F'({\mathcal
O}(k))\stackrel{\sim}{\lto} \Phi_{E'}({\mathcal O}(k)),$ completing
the following commutative diagram
$$
\begin{array}{ccc}
{\mr H}^0 ({\mathcal O}(k))\ot F'({\mathcal
O})&\stackrel{F'(\delta_k)}{\lto}
&F'({\mathcal O}(k))\\
\llap{\ss{id}}\da&&\da\rlap{\ss{f_k}}\\
{\mr H}^0 ({\mathcal O}(k))\ot F'({\mathcal O})&\stackrel{\pi'_*
(\gamma_k)}{\lto}&\Phi_{E'}({\mathcal O}(k))
\end{array}
$$

Now we have to show that these morphisms are functorial. For any
$\alpha : {\mathcal O}(k)\to{\mathcal O}(l)$ we have the commutative
squares
$$
\begin{array}{ccc}
{\mr H}^0 ({\mathcal O}(k))\ot F'({\mathcal
O})&\stackrel{F'(\delta_k)}{\lto}
&F'({\mathcal O}(k))\\
\llap{\ss{{\mr H}^0 (\alpha)\ot id}}\da&&\da\rlap{\ss{F'(\alpha)}}\\
{\mr H}^0 ({\mathcal O}(l))\ot F'({\mathcal
O})&\stackrel{F'(\delta_l)}{\lto}& F'({\mathcal O}(l))
\end{array}
$$
and
$$
\begin{array}{ccc}
{\mr H}^0 ({\mathcal O}(k))\ot F'({\mathcal O})& \stackrel{\pi'_*
(\gamma_k)}{\lto}
&\Phi_{E'}({\mathcal O}(k))\\
\llap{\ss{{\mr H}^0 (\alpha)\ot id}}\da&&
\da\rlap{\ss{\Phi_{E'}(\alpha)}}\\
{\mr H}^0 ({\mathcal O}(l))\ot F'({\mathcal O})& \stackrel{\pi'_*
(\gamma_l)}{\lto}&\Phi_{E'}({\mathcal O}(l))
\end{array}
$$
Therefore we have the equalities:
$$
\begin{array}{l}
f_l\circ F'(\alpha)\circ F'(\delta_k)=f_l\circ F'(\delta_l)\circ
 ({\mr H}^0 (\alpha)\ot id)=
\pi'_* (\gamma_l)\circ
 ({\mr H}^0 (\alpha)\ot id)=
\Phi_{E'}(\alpha)\circ \pi'_* (\gamma_k)= \\
\Phi_{E'}(\alpha)\circ f_k \circ F'(\delta_k)
\end{array}
$$
Since the complexes $D^{\bdot}_k$ and $D^{\bdot}_l$ satisfy the
conditions of Lemma \ref{pd2} there exists only one morphism $h:
F'({\mathcal O}(k)) \to \Phi_{E'}({\mathcal O}(l))$ such that
$$
h\circ F'(\delta_k) = \pi'_* (\gamma_l)\circ
 ({\mr H}^0 (\alpha)\ot id)
$$
Hence $f_l\circ F'(\alpha)$ coincides with $\Phi_{E'}(\alpha)\circ f_k.$

Now, consider the case $k<0.$

Let us take the following right resolution for ${\mathcal O}(k)$ on
${\Bbb{P}}^N.$
$$
{\mathcal O}(k)\stackrel{\sim}{\lto}\bigl\{ V^k_0 \ot{\mathcal
O}\lto\cdots\lto V^k_N \ot {\mathcal O}(N)\bigr\}
$$
By Lemma \ref{pd2}, the morphism of the complexes over $\db{X}$
$$
\begin{array}{ccccc}
V^k_0 \ot F'({\mathcal O})&\lto&\cdots&\lto& V^k_N \ot F'({\mathcal O}(N))\\
\llap{\ss{id\ot f_0}}\da\wr&&&&\llap{\ss{id\ot f_N}}\da\wr\\
V^k_0 \ot \Phi_{E'}({\mathcal O})&\lto&\cdots&\lto& V^k_N \ot \Phi_{E'}({\mathcal O}(N))\\
\end{array}
$$
gives us the uniquely determined morphism $f_k : F'({\mathcal
O}(k))\lto \Phi_{E'}({\mathcal O}(k)).$

It is not hard to prove that these morphisms are functorial. The proof
 is left to a reader. $\Box$

\raz Now we must prove that there exists an  object $E\in\db{M\times X}$
such that $j_* E\cong E'.$

Let ${\mathcal L}$ be a very ample invertible sheaf on $M$ and let
$j: M\hookrightarrow {\Bbb{P}}^N$ be an embedding with respect to
${\mathcal L}.$ By $A$ denote the graded algebra
$\bigoplus\limits^{\infty}_{i=0} {\mr H}^0 (M, {\mathcal L}^{i}).$

Let $B_0=k,$ and $B_1=A_1.$ For $m\ge 2,$ we define $B_m$ as
\begin{equation} \label{dual}
B_m = Ker( B_{m-1}\ot A_1 \lto B_{m-2}\ot A_2 )
\end{equation}
\th{Definition}
$A$ is said to be $n$-Koszul if the following sequence is exact
$$
B_n \ot_k A\lto B_{n-1}\ot_k A\lto \cdots\lto B_1\ot_k A\lto A\lto k\lto 0
$$
\eth

Assume that $A$ is n-Koszul. Let $R_0 = {\mathcal O}_M.$ For $m\ge
1,$ denote by $R_m$ the kernel of the morphism $B_m\ot{\mathcal O}_M
\lto B_{m-1}\ot{\mathcal L}.$ Using (\ref{dual}), we obtain the
canonical morphism $R_m\lto A_1\ot R_{m-1}.$ (actually, ${\h R_m,
{R_{m-1}}}\cong A_1^*$).

Since $A$ is $n$-Koszul, we have the exact sequences
$$
0\lto R_m \lto B_m\ot {\mathcal O}_M \lto B_{m-1}\ot{\mathcal
L}\lto\cdots \lto B_1 \ot {\mathcal L}^{{m-1}}\lto {\mathcal
L}^{m}\lto 0
$$
for $m\le n.$

We have the canonical morphisms $f_m : j^* \Omega^m (m)\lto R_m,$ because
$\Lambda^i A_1 \subset B_i$ and there exist the exact sequences on ${\Bbb{P}}^N$
$$
0\lto \Omega^m (m) \lto \Lambda^m A_1 \ot {\mathcal O}\lto
\Lambda^{m-1} A_1 \ot {\mathcal O}(1)\lto \cdots \lto {\mathcal
O}(m)\lto 0
$$

It is known that for any $n$ there exists $l$ such that the Veronese
algebra $A^l = \bigoplus\limits^{\infty}_{i=0} {\mr H}^0 (M,
{\mathcal L}^{il})$ is $n$-Koszul.( Moreover, it was proved in
\cite{Bec} that $A^l$ is Koszul for $l\gg 0$).

 Using the technique of \cite{IM} and substituting ${\mathcal L}$ with ${\mathcal L}^j,$ when $j$ is sufficiently large ,  we can choose for any $n$ a very ample ${\mathcal L}$ such that

1) algebra $A$
is $n$-Koszul,

2) the  complex
$$
{\mathcal L}^{ -n}\bt R_n\lto\cdots\lto {\mathcal L}^{ -1}\bt
R_1\lto {\mathcal O}_M \bt R_0 \lto {\mathcal O}_{\Delta}
$$
 on $M\times M$ is exact,

3) the following
sequences on $M.$
$$
A_{k-n}\ot R_n \lto A_{k-n+1}\ot R_{n-1}\lto\cdots\lto A_{k-1}\ot
R_1 \lto A_k \ot R_0 \lto {\mathcal L}^{k}\lto 0
$$
are exact  for any $k\ge 0.$ Here, by definition, if $ k-i<0,$ then $A_{k-i}=0.$
(see Appendix for proof).

Let us denote by $T_k$ the kernel
of the morphism $A_{k-n}\ot R_n \lto A_{k-n+1}\ot R_{n-1}.$

Consider the following complex over $\db{M\times X}$
\begin{equation}\label{nob}
{\mathcal L}^{ -n}\bt F(R_n)\lto\cdots\lto {\mathcal L}^{ -1}\bt
F(R_1)\lto {\mathcal O}_M \bt F(R_0)
\end{equation}
Here the morphism ${\mathcal L}^{-k}\bt F(R_k)\lto {\mathcal
L}^{-k+1}\bt F(R_{k-1})$ is induced by the canonical morphism
$R_k\lto A_1 \ot R_{k-1}$ with respect to the following sequence of
isomorphisms
$$
{\h {\mathcal L}^{-k}\bt F(R_k), {{\mathcal L}^{-k+1}\bt
F(R_{k-1})}}\cong {\h  F(R_k), {{\mr H}^0 ({\mathcal L})\ot
F(R_{k-1})}}\cong
$$
$$
\cong{\h  R_k, {A_1 \ot R_{k-1}}}
$$

By Lemma \ref{pd1}, there is a convolution of the complex (\ref{nob}) and
all convolutions are isomorphic. Let $G\in \db{M\times X}$ be a convolution
of this complex.

For any $k\ge 0,$ object $\pi_* (G\ot p^* ({\mathcal L}^k))$ is a
convolution of the complex
$$
A_{k-n}\ot F(R_n)\lto A_{k-n+1}\ot F(R_{n-1})\lto\cdots\lto A_k \ot F(R_0).
$$
On the other side, we know that $T_k[n]\oplus{\mathcal L}^k$ is a
convolution of the complex
$$
A_{k-n}\ot R_n\lto A_{k-n+1}\ot R_{n-1}\lto\cdots\lto A_k \ot R_0,
$$
because ${\E n+1 , {}, {\mathcal L}^k , {T_k}}=0$ for $n\gg 0.$
Therefore, by Lemma \ref{pd1}, we have $\pi_* (G\ot p^* ({\mathcal
L}^k))\cong F(T_k[n]\oplus {\mathcal L}^k).$

It follows immediately from Remark \ref{boun} that  the cohomology
sheaves ${H}^i (\pi_* (G\ot p^* ({\mathcal L}^k)))= {H}^i
(F(T_k)[n])\oplus{H}^i (F( {\mathcal L}^k))$ concentrate on the
union $[-n-a, -n]\cup [-a, 0]$ for any $k>0$ ($a$ was defined in
\ref{boun}). Therefore the cohomology sheaves ${H}^i (G)$
 also concentrate on $[-n-a, -n]\cup [-a, 0].$ We can assume that $n> dimM
+ dimX + a.$ This implies that $G\cong C\oplus E,$ where $E, C$ are
objects of $\db{M\times X}$ such that ${H}^i (E)=0$ for $i\not\in
[-a, 0]$ and ${H}^i (C)=0$ for $i\not\in [-n-a, -n].$ Moreover, we
have $\pi_* (E\ot p^* ({\mathcal L}^k))\cong F({\mathcal L}^k).$

Now we show that $j_* (E)\cong E'.$
Let us consider the morphism of the complexes over $\db{{\Bbb{P}}^N \times X}.$
$$
\begin{array}{ccccc}
{\mathcal O}(-n)\bt F' (\Omega^n (n))&\lto&\cdots&\lto&{\mathcal O}\bt F' ({\mathcal O})\\
\da\rlap{\ss{can\bt F(f_n)}}&&&&\da\rlap{\ss{can\bt F(f_0)}}\\
j_* ({\mathcal L}^{-n}) \bt F(R_n)&\lto&\cdots&\lto& j_* ({\mathcal
O}_M)\bt F(R_0)
\end{array}
$$
By Lemma \ref{pd2}, there exists a morphism of  convolutions $\phi :
K\lto j_* (G).$ If $N>n,$ then $K$ is not isomorphic to $E',$ but there is
a distinguished triangle
$$
S\lto K\lto E'\lto S[1]
$$
and the cohomology sheaves ${H}^i (S)\ne 0$ only if $i\in [-n-a, -n].$
Now, since ${\h S, {j_* (E)}}=0$ and ${\h S[1], {j_* (E)}}=0,$ we have a uniquely
determined morphism $\psi : E' \lto j_* (E)$ such that the following diagram
commutes
$$
\begin{array}{ccc}
K&\stackrel{\phi}{\lto}& j_* (G)\\
\da&&\da\\
E'&\stackrel{\psi}{\lto}& j_* (E)
\end{array}
$$
We know that $\pi'_*( E'\ot q^*({\mathcal O}(k)))\cong F({\mathcal
L}^k)\cong
 \pi_*( E\ot p^*({\mathcal L}^k)).$ Let $\psi_k$ be the morphism $
\pi'_*( E'\ot q^*({\mathcal O}(k)))\lto \pi_*( E\ot p^*({\mathcal
L}^k))$ induced by $\psi.$ The morphism $\psi_k$ can be included in
the following commutative diagram:
$$
\begin{array}{ccccc}
S^k A_1 \ot F({\mathcal O})&\stackrel{can}{\lto}& F({\mathcal
L}^k)&\stackrel{\sim}{\lto}
&\pi'_*( E'\ot q^*({\mathcal O}(k)))\\
\llap{\ss{can}}\da&&&&\da\rlap{\ss{\psi_k}}\\
A_k \ot F({\mathcal O})&\stackrel{can}{\lto}& F({\mathcal
L}^k)&\stackrel{\sim}{\lto} &\pi_*( E\ot p^*({\mathcal L}^k))
\end{array}
$$
Thus we  see that $\psi_k$ is an isomorphism for any $k\ge 0.$ Hence
$\psi$ is an isomorphism too. This proves the following:
\th{Lemma}\label{obj} There exists an object $E\in\db{M\times X}$ such that
$j_* (E)\cong E',$ where $E'$ is the object from $\db{{\Bbb{P}}^N \times X},$
constructed in \ref{cons}.
\eth

%%%%%%%%%%%%%%%%%%%%%%%%%%%%%%%%%%%%%%%%%%%%%%%%%%%%%%%%%%%%%%%%%%%%%

\raz
Now, we prove some statements relating to abelian categories.
they are needed for the sequel.

Let  ${\mathcal A}$ be a $k$-linear  abelian category (henceforth we
shall consider only $k$-linear abelian categories). Let $\{ P_i
\}_{i\in\Bbb Z}$ be a sequence of objects from ${\mathcal A}.$
\th{Definition} We say that this sequence is  {\sf ample} if for
every object $X\in {\mathcal A}$ there exists $N$  such that for all
$i<N$ the following conditions hold:

a) the canonical morphism ${\h P_i, X}\otimes P_i \lto X$ is surjective,

b) ${\E j, {}, P_i, X}=0$ for any $j\not=0,$

c) ${\h X, {P_i}}=0.$ \eth It is clear that if ${\mathcal L}$ is  an
ample invertible sheaf on a projective variety in usual sense, then
the
 sequence $\{ {\mathcal L}^i \}_{i\in{\Bbb{Z}}}$ in the abelian
 category of coherent sheaves  is ample.

\th{Lemma}\label{zer1} Let $\{ P_i \}$ be an ample sequence in an
abelian category ${\mathcal A}.$ If $X$ is an object in
$\db{\mathcal A}$ such that ${\H {\bdot}, {}, P_i, X}=0$ for all
$i\ll 0,$ then $X$ is the zero object. \eth \pr If $i\ll 0$ then
$${\h P_i, {{H}^k (X)}}\cong{\H k, {}, P_i, X}=0$$
The morphism ${\h P_i, {{H}^k (X)}}\otimes P_i\lto {H}^k (X)$ must be
surjective for $i\ll 0,$ hence  ${H}^k (X)=0$ for all $k.$ Thus
$X$ is the zero object. $\Box$

\th{Lemma}\label{zer2} Let $\{ P_i \}$ be an ample sequence in an
abelian category ${\mathcal A}$ of finite homological dimension. If
$X$ is an object in $\db{\mathcal A}$ such that ${\H {\bdot}, {}, X,
{P_i}}=0$ for all $i\ll 0.$ Then $X$ is the zero object. \eth \pr
Assume that the cohomology objects of $X$ are concentrated in a
segment $[a, 0].$ There exists the canonical morphism $X\lto {H}^0
(X).$ Consider a surjective morphism $P_{i_1}^{\oplus k_1}\lto {H}^0
(X).$ By $Y_1$ denote the kernel of this morphism. Since ${\H
{\bdot}, {}, X, {P_{i_1}}}=0$ we have ${\H 1, {}, X, {Y_1}}\not=0.$
Further take a surjective morphism $P_{i_2}^{\oplus k_2}\lto Y_1.$
By $Y_2$ denote the kernel of this morphism. Again,  since ${\H
{\bdot}, {}, X, {P_{i_2}}}=0,$ we obtain ${\H 2, {}, X,
{Y_2}}\not=0.$  Iterating this procedure as needed, we get
contradiction with the assumption that ${\mathcal A}$ is of finite
homological dimension. $\Box$

\th{Lemma}\label{f&f} Let ${\mathcal B}$ be an abelian category,
${\mathcal A}$ an abelian category of finite homological dimension,
and $\{ P_i \}$ an ample sequence in ${\mathcal A}.$ Suppose $F$ is
an exact functor from $\db{\mathcal A}$ to $\db{\mathcal B}$ such
that it has right and left adjoint functors $F^!$ and $F^*$
respectively. If the maps
$$
{\H k, {}, P_i, {P_j}}\stackrel{\sim}{\lto}{\H k, {}, F(P_i), {F(P_j)}}
$$
are isomorphisms for $i<j$ and all $k.$
Then $F$ is full and faithful.
\eth
\pr
Let us take the canonical morphism $f_i : P_i\lto F^! F(P_i)$ and consider
a distinguished triangle
$$
P_i\stackrel{f_i}{\lto} F^! F(P_i)\lto C_i\lto P_i [1].
$$
Since for $j\ll 0$ we have isomorphisms:
$$
{\H k, {}, P_j, {P_i}}\stackrel{\sim}{\lto}{\H k, {}, F(P_j), {F(P_i)}}
\cong{\H k, {}, P_j, {F^! F(P_i)}}.
$$
We see that ${\H {\bdot}, {}, P_j, {C_i}}=0$ for $j\ll 0.$
It follows from Lemma \ref{zer1} that $C_i =0.$ Hence $f_i$ is
an isomorphism.

Now, take the canonical morphism $g_X : F^* F(X)\lto X$ and consider
a distinguished triangle
$$
F^* F(X)\stackrel{g_X}{\lto} X\lto C_X\lto F^* F(X)[1]
$$
We have the following sequence of isomorphisms
$$
{\H k, {}, X, {P_i}}\stackrel{\sim}{\lto}{\H k, {}, X, {F^! F(P_i)}}
\cong{\H k, {}, F^* F(X), {P_i}}
$$
This implies that ${\H {\bdot}, {}, C_X, {P_i}}=0$ for all $i.$
By Lemma \ref{zer2}, we obtain $C_X=0.$ Hence $g_X$ is an isomorphism.
It follows that $F$ is full and faithful. $\Box$

Let ${\mathcal A}$ be an abelian category possessing an ample
sequence $\{ P_i\}.$ Denote by $\db{\mathcal A}$  the bounded
derived category of ${\mathcal A}.$ Let us consider the full
subcategory $j: {\mathcal C}\hookrightarrow\db{\mathcal A}$ such
that $\mbox{Ob} {\mathcal C}:=\{  P_i\; |\; i\in{\Bbb{Z}} \}.$ Now
we would like to show that if there exists an isomorphism of a
functor $F : \db{{\mathcal A}}\lto\db{{\mathcal A}}$ to identity
functor on the subcategory ${\mathcal C},$ then it can be extended
to the whole $\db{\mathcal A}.$

\th{Proposition}\label{ext} Let $F :\db{{\mathcal
A}}\lto\db{{\mathcal A}}$ be an autoequivalence. Suppose there
exists an isomorphism $f : j\stackrel{\sim}{\lto}F\mid_{\mathcal C}
$ ( where $j : {\mathcal C}\hookrightarrow \db{{\mathcal A}}$ is  a
natural embedding).
 Then it can be extended to an isomorphism
$id\stackrel{\sim}{\lto}F$ on the whole $\db{\mathcal A}.$ \eth \pr
First, we can extend the transformation $f$ to all direct sums of
objects ${\mathcal C}$ componentwise , because $F$ takes direct sums
to direct sums.

Note that $X\in\db{\mathcal A}$ is isomorphic to an object in
${\mathcal A}$ iff ${\H j, {}, P_i, X}=0$ for $j\not=0$ and $i\ll
0.$ It follows that $F(X)$ is isomorphic to an object in ${\mathcal
A},$ because
$$
{\H j, {}, P_i, {F(X)}}\cong{\H j, {}, F(P_i), {F(X)}}\cong
{\H j, {}, P_i, X}=0
$$
for $j\not=0$ and $i\ll 0.$

\praz At first, let $X$ be an object from ${\mathcal A}.$ Take a
surjective morphism $v : P^{\oplus k}_i\lto X.$ We have the morphism
$f_i : P^{\oplus k}_i\lto F(P^{\oplus k}_i)$ and two distinguished
triangles:
$$
\begin{array}{ccccccc}
Y&\stackrel{u}{\lto}&P^{\oplus k}_i&\stackrel{v}{\lto}&X&\lto&Y [1]\\
&&\da\rlap{$f_i$}&&&& \\
F(Y)&\stackrel{F(u)}{\lto}&F(P^{\oplus k}_i)&\stackrel{F(v)}{\lto}&F(X)&\lto&F(Y)[1]\\
\end{array}
$$

Now we show that $F(v)\circ f_i\circ u=0.$
Consider any surjective morphism $ w : P^{\oplus l}_j\lto Y.$
It is sufficient to check that $F(v)\circ f_i\circ u\circ w=0.$
Let $f_j : P^{\oplus l}_j\lto F(P^{\oplus l}_j)$ be the canonical morphism.
Using the commutation relations for $f_i$ and $f_j,$ we obtain
$$
F(v)\circ f_i\circ u\circ w = F(v)\circ F(u\circ w)\circ f_j=
F(v\circ u\circ w)\circ f_j =0
$$
because $v\circ u=0.$

Since ${\h Y[1], {F(X)}}=0,$ by Lemma \ref{tr}, there exists a unique morphism
$f_X : X\lto F(X)$ that commutes with $f_i.$

\praz Let us show that $f_X$ does not depend from morphism $v : P^{\oplus k}_i
\lto X.$
Consider two surjective morphisms $v_1 : P^{\oplus k_1}_{i_1}\lto X$ and
$v_2 : P^{\oplus k_2}_{i_2}\lto X.$ We can take two surjective
morphisms $w_1 : P^{\oplus l}_{j}\lto P^{\oplus k_1}_{i_1}$ and
 $w_2 : P^{\oplus l}_{j}\lto P^{\oplus k_2}_{i_2}$ such that the
following diagram is commutative:
$$
\begin{array}{ccc}
P^{\oplus l}_{j}&\stackrel{w_2}{\lto}& P^{\oplus k_2}_{i_2}\\
 \da\rlap{\ss{w_1}}&& \da\rlap{\ss{v_2}}\\
P^{\oplus k_1}_{i_1}&\stackrel{v_1}{\lto}& X\\
\end{array}
$$
Clearly, it is sufficient to check the coincidence of the morphisms,
constructed by $v_1$ and $v_1\circ w_1.$ Now, let us consider the following
commutative diagram:
$$
\begin{array}{ccccc}
 P^{\oplus l}_{j}&\stackrel{w_1}{\lto}& P^{\oplus k_1}_{i_1}&
\stackrel{v_1}{\lto}&X\\
\da\rlap{\ss{f_j}}&&\da\rlap{\ss{v_2}}&&\da\rlap{\ss{f_X}}\\
F(P^{\oplus l}_{j})&\stackrel{F(w_1)}{\lto}&F(P^{\oplus k_1}_{i_1})&
\stackrel{F(v_1)}{\lto}&F(X)\\
\end{array}
$$
Here the morphism $f_X$ is constructed by $v_1.$ Both squares of this
diagram are
commutative. Since there exists only one morphism from $X$ to $F(X)$ that
commutes with $f_j,$ we see that the $f_X,$ constructed by $v_1,$ coincides
with the morphism, constructed by
$v_1\circ w_1.$

\praz Now we must show that for any morphism $X\stackrel{\phi}{\lto}Y$
we have equality:
$$
f_Y \circ \phi = F(\phi )\circ f_X
$$
Take a surjective morphism  $P^{\oplus l}_j\stackrel{v}{\lto} Y.$
 Choose a surjective morphism $P^{\oplus k}_i\stackrel{u}{\lto}X$
 such that the composition $\phi\circ u$  lifts to a morphism
$\psi : P^{\oplus k}_i{\lto} P^{\oplus l}_j .$ We can do it, because
for $i\ll 0$ the map ${\h P^{\oplus k}_i, {P^{\oplus l}_j}}\to
{\h P^{\oplus k}_i, Y}$ is surjective. We get the commutative square:
$$
\begin{array}{ccc}
P^{\oplus k}_{i}&\stackrel{u}{\lto}&X\\
\da\rlap{\ss{\psi}}&&\da\rlap{\ss{\phi}}\\
P^{\oplus l}_{j}&\stackrel{v}{\lto}&Y\\
\end{array}
$$
By $h_1$ and $h_2$ denote $f_Y\circ\phi$ and $F(\phi)\circ f_X$ respectively.
We have the following sequence of equalities:
$$
h_1\circ u=f_Y\circ\phi\circ u=f_Y\circ v\circ\psi=F(v)\circ f_j\circ\psi=
F(v)\circ F(\psi)\circ f_i
$$
and
$$
h_2\circ u=F(\phi)\circ f_X\circ u =F(\phi)\circ F(u)\circ f_i =
F(\phi\circ u)\circ f_i = F(v\circ \psi)\circ f_i = F(v)\circ F(\psi)\circ f_i
$$
Consequently, the following square is commutative for $t=1,2.$
$$
\begin{array}{ccccccc}
Z&\lto&P^{\oplus k}_{i}&\stackrel{u}{\lto}&X&\lto&Z[1]\\
&&\llap{\ss{F(\psi)\circ f_i}}\da&&\da\rlap{\ss{h_t}}&&\\
F(W)&\lto&F(P^{\oplus l}_{j})&\stackrel{F(v)}{\lto}&F(Y)&\lto&F(W)[1]\\
\end{array}
$$
By Lemma \ref{tr}, as ${\h Z[1], {F(Y)}}=0,$ we obtain  $h_1=h_2.$
Thus, $f_Y \circ \phi = F(\phi )\circ f_X.$

Now take a cone $C_X$ of the morphism $f_X.$ Using the following
isomorphisms
$$
{\h P_i, X}\cong{\h F(P_i), {F(X)}}\cong{\h P_i, {F(X)}},
$$
we obtain ${\H j, {}, P_i, {C_X}}=0$ for all $j,$ when $i\ll 0.$
Hence, by Lemma \ref{zer1}, $C_X=0$ and $f_X$ is an isomorphism.

%%%%%%%%%%%%%%%%%%%%%%%%%%%%%%%%%%%%%%%%%%%%%%%%%%%%%%%
%%%%%%%%%%%%%%%%%%%%%%%%%%%%%%%%%%%%%%%%%%%%%%%%%%%%%%%
\praz Let us define $f_{X[n]} : X[n]\lto F(X[n])\cong F(X)[n]$ for
any $X\in{\mathcal A}$ by
$$
f_{X[n]}=f_X [n].
$$

It is easily shown that these transformations commute with any
$u\in{\E k, {}, X, Y}.$
Indeed, since any element $u\in{\E k, {}, X, Y}$ can be represented as a composition $ u = u_0 u_1 \cdots u_k $
of some elements  $u_i\in{\E 1, {}, Z_i, {Z_{i+1}}}$ and $Z_0=X, Z_k=Y,$
we have only to check it for $u\in{\E 1, {}, X, Y}$ .
Consider  the following diagram:
$$
\begin{array}{ccccccc}
Y&\lto&Z&\lto&X&\stackrel{u}{\lto}&Y[1]\\
\llap{\ss{f_Y}}\da&&\da\rlap{\ss{f_Z}}&&&&\da\rlap{\ss{f_Y [1]}}\\
F(Y)&\lto&F(Z)&\lto&F(X)&\stackrel{F(u)}{\lto}&F(Y)[1]\\
\end{array}
$$
By an axiom of triangulated categories there exists a morphism $h : X\to F(X)$
such that $(f_Y, f_Z, h)$ is a morphism of triangles.
On the other hand, since ${\h Y[1], {F(X)}}=0,$ by Lemma \ref{tr}, $h$ is a unique morphism
that commutes with $f_Z.$
But $f_X$ also commutes with $f_Z.$ Hence we have $h=f_X.$ This implies that
$$
f_Y [1] \circ u = F(u) \circ f_X
$$

\praz The rest of  the proof is by induction over the length of a
segment, in which the cohomology objects of $X$ are concentrated.
Let $X$ be  an object from $ \db{\mathcal A}$ and suppose that its
cohomology objects ${H}^p (X)$ are concentrated in a segment  $[a,
0].$ Take $v: P^{\oplus k}_i\lto X$ such that
\begin{eqnarray}\label{tt}
a)& {\H j, {}, P_i, {{H}^p (X)}}=0& \mbox{ for all }\; p\; \mbox{ and for }
\; j\not=0,\nonumber\\
b)& u: P_i^{\oplus k}\lto {H}^0 (X)& \mbox{is  the surjective morphism},\\
c)&{\h {H}^0 (X), {P_i}}=0.&\nonumber
\end{eqnarray}
Here $u$ is the composition
$v$ with the canonical morphism $X\lto {H}^0 (X).$
Consider a distinguished triangle:
$$
Y[-1]\lto P^{\oplus k}_i\stackrel{v}{\lto}X\lto Y
$$
By the induction hypothesis, there exists the isomorphism $f_Y$ and it commutes
with $f_i.$ So we have the commutative diagram:
$$
\begin{array}{ccccccc}
Y[-1]&\lto&P^{\oplus k}_i&\stackrel{v}{\lto}&X&\lto&Y\\
\llap{\ss{f_Y [-1]}}\da&&\da\rlap{\ss{f_i}}&&&&\da\rlap{\ss{f_Y }}\\
F(Y)[-1]&\lto&F(P^{\oplus k}_i)&\stackrel{F(v)}{\lto}&F(X)&\lto&F(Y)\\
\end{array}
$$
Moreover we have the following sequence of equalities
$$
{\h X, {F(P^{\oplus k}_i)}}\cong{\h X, {P^{\oplus k}_i}}\cong
{\h {H}^0 (X), {P^{\oplus k}_i}}=0
$$
Hence, by Lemma \ref{tr}, there exists a unique morphism $f_X : X\lto F(X)$ that commutes with $f_Y.$

\praz  We must first show that $f_X$ is correctly defined.
Suppose we have two morphisms $v_1 : P_{i_1}^{\oplus k_1}\lto X$
and $v_2 : P_{i_2}^{\oplus k_2}\lto X.$ As above, we can find $P_j$ and
surjective morphisms $w_1, w_2$ such that the following diagram is commutative
$$
\begin{array}{ccc}
P^{\oplus l}_{j}&\stackrel{w_2}{\lto}& P^{\oplus k_2}_{i_2}\\
 \da\rlap{\ss{w_1}}&& \da\rlap{\ss{u_2}}\\
P^{\oplus k_1}_{i_1}&\stackrel{u_1}{\lto}&{H}^0 (X)\\
\end{array}
$$
We can find a morphism $\phi : Y_j\lto Y_{i_1}$ such that  the triple
$(w_1, id, \phi)$ is a morphism of distinguished triangles.
$$
\begin{array}{ccccccc}
P^{\oplus l}_j&\stackrel{v_1\circ w_1}{\lto}&X&\lto&Y_j&\lto&P^{\oplus l}_j [1]\\
\llap{\ss{w_1}}\da&&\da\rlap{\ss{id}}&&\da\rlap{\ss{\phi}}&&\da\rlap{\ss{w_1 [1]}}\\
P^{\oplus k_1}_{i_1}&\stackrel{v_1}{\lto}&X&\lto&Y_{i_1}&\lto&
 P^{\oplus k_1}_{i_1} [1]\\
\end{array}
$$
By the induction hypothesis, the following square is commutative.
$$
\begin{array}{ccc}
Y_j&\stackrel{\phi}{\lto}&Y_{i_1}\\
\llap{\ss{f_{Y_j}}}\da&&\da\rlap{\ss{f_{Y_i}}}\\
F(Y_j)&\stackrel{F(\phi)}{\lto}&F(Y_{i_1})\\
\end{array}
$$
Hence, we see that the $f_X,$ constructed by $v_1\circ w_1,$ commutes with
$f_{Y_{i_1}}$ and, consequently, coincides with the $f_X,$ constructed by $v_1$;
because such morphism is unique by Lemma \ref{tr}. Therefore morphism $f_X$ does not
depend on a choice of morphism $v: P^{\oplus k}_i\lto X.$

\praz Finally, let us prove that for any morphism $\phi : X\lto Y$ the following
 diagram commutes
\begin{equation}\label{comd}
\begin{array}{ccc}
X&\stackrel{\phi}{\lto}&Y\\
\llap{\ss{f_X}}\da&&\da\rlap{\ss{f_Y }}\\
F(X)&\stackrel{F(\phi)}{\lto}&F(Y)
\end{array}
\end{equation}
Suppose the cohomology objects of $X$ are concentrated on a segment $[a, 0]$
and the cohomology objects of $Y$ are concentrated on $[b, c].$

\noindent{\it Case 1.} If $c<0,$ we take a morphism $v: P^{\oplus k}_i
\lto X$  that satisfies conditions (\ref{tt}) and ${\h P^{\oplus k}_i, Y}=0.$
We have a distinguished triangle:
$$
\begin{array}{ccccccc}
P^{\oplus k}_i&\stackrel{v_1}{\lto}&X&\stackrel{\alpha}{\lto}&Z&\lto&
P^{\oplus k}_i [1]\\
\end{array}
$$
Applying the functor ${\h -, Y}$ to this triangle we found that there exists a morphism $\psi : Z \lto Y$ such that $\phi = \psi \circ \alpha.$
We know that $f_X,$ defined above, satisfy
$$
F(\alpha )\circ f_X = f_Z \circ \alpha
$$
If we assume that the diagram
$$
\begin{array}{ccc}
Z&\stackrel{\psi}{\lto}&Y\\
\llap{\ss{f_Z}}\da&&\da\rlap{\ss{f_Y }}\\
F(Z)&\stackrel{F(\psi)}{\lto}&F(Y)\\
\end{array}
$$
commutes, then  diagram (\ref{comd}) commutes too.

This means that for verifying the commutativity of (\ref{comd}) we can
substitute $X$ by an object $Z.$
And the cohomology objects of $Z$ are concentrated on the segment $[a, -1].$

\noindent{\it Case 2.} If $c\ge0,$ we take a surjective morphism
$v: P^{\oplus k}_i\lto Y[c]$ that satisfies conditions (\ref{tt}) and
 ${\h {H}^c (X), {P^{\oplus k}_i}}=0.$
Consider a distinguished triangle
$$
\begin{array}{ccccccc}
P^{\oplus k}_i [-c]&\stackrel{v[-c]}{\lto}&Y&\stackrel{\beta}{\lto}&W&\lto&
P^{\oplus k}_i [-c+1]\\
\end{array}
$$
Note that the cohomology objects of $W$ are concentrated on $[b, c-1].$

By $\psi$ denote the composition $\beta\circ\phi.$ If we assume that the following
square
$$
\begin{array}{ccc}
X&\stackrel{\psi}{\lto}&W\\
\llap{\ss{f_X}}\da&&\da\rlap{\ss{f_W }}\\
F(X)&\stackrel{F(\psi)}{\lto}&F(W)\\
\end{array}
$$
commutes, then, since $F(\beta)\circ f_Y = f_W \circ \beta,$
 $$
F(\beta)\circ(f_Y\circ\phi - F(\phi)\circ f_X)= f_W \circ \psi -   F(\psi ) \circ f_X =0.
$$
We chose $P_i$ such that ${\h X, {P^{\oplus k}_i [-c]}}=0.$ As $F( P_i^{\oplus k})$ is isomorphic to $P_i^{\oplus k},$ then ${\h X, {F(P^{\oplus k}_i [-c])}}=0.$
Applying the functor ${\h X, {F(-)}}$ to the above triangle we found that the composition with $F(\beta )$ gives an inclusion of ${\h X, {F(Y)}}$ into ${\h X, {F(W)}}.$ This follows that $f_Y\circ\phi = F(\phi)\circ f_X,$ i.e. our diagram (\ref{comd}) commutes.

Combining case 1 and case 2, we can reduce the checking of
commutativity of diagram (\ref{comd}) to the case when $X$ and $Y$
are objects from the abelian category ${\mathcal A}.$ But for those
it has already been done. Thus the proposition is proved. $\Box$

\raz {\bf Proof of theorem}. 1) {\sc Existence}. Using Lemma
\ref{obj} and Lemma \ref{tran}, we can construct an object
$E\in\db{M\times X}$ such that there exists an isomorphism of the
functors $\bar{f} : F\bigl|_{\mathcal C}\stackrel{\sim}{\lto}\Phi_E
\bigr|_{\mathcal C}$ on full subcategory ${\mathcal
C}\subset\db{M},$ where ${\mr O}{\mr b}{\mathcal C}= \{ {\mathcal
L}^{i} \mid i\in{\Bbb{Z}} \}$ and ${\mathcal L}$ is a very ample
invertible sheaf on $M$ such that for any $k>0$ ${\rm H}^i ( M,
{\mathcal L}^k )=0$ , when $i\ne 0.$

By Lemma \ref{f&f} the functor $\Phi_E$ is full and faithfull. Moreover, the functors $F^!\circ \Phi_E$ and
$\Phi_E^* \circ F$ are full and
faithful too, because we have the isomorphisms:
$$
F^! (\bar{f}) : F^!\circ F\bigl|_{\mathcal C}\cong id_{\mathcal
C}\stackrel{\sim}{\lto}F^! \circ \Phi_E\bigl|_{\mathcal C}
$$
$$
\Phi_E^* (\bar{f}) : \Phi_E^*\circ F\bigl|_{\mathcal
C}\stackrel{\sim}{\lto} \Phi_E^*\circ \Phi_E\bigl|_{\mathcal C}\cong
id_{\mathcal C}
$$
and conditions of  Lemma \ref{f&f} is fulfilled.

Further, the functors $F^!\circ \Phi_E$  and
$\Phi_E^* \circ F$ are  equivalences,
because they are adjoint each other.

Consider the isomorphism $F^! (\bar{f}) : F^!\circ F\bigl|_{\mathcal
C}\cong id_{\mathcal C}\stackrel{\sim}{\lto}F^! \circ
\Phi_E\bigl|_{\mathcal C}$ on the subcategory ${\mathcal C}.$ By
Proposition \ref{ext} we can extend it onto the whole $\db{M},$ so
$id\stackrel{\sim}{\lto}F^!\circ\Phi_E.$

Since $F^!$ is the right adjoint to $F,$ we get the morphism of the
functors $f : F\lto \Phi_E$ such that $f|_{\mathcal C} =\bar{f}.$ It
can easily be checked that $f$ is an isomorphism. Indeed, let $C_Z$
be a cone of the morphism $f_Z : F(Z)\lto \Phi_E (Z).$ Since $F^!
(f_Z)$ is an isomorphism, we obtain $F^! (Z)=0.$ Therefore, this
implies that ${\h F(Y), {C_Z}}=0$ for any object $Y.$ Further, there
are isomorphisms $F({\mathcal L}^k)\cong \Phi_E ({\mathcal L}^k)$
for any $k.$ Hence, we have
$$
{\H i, {},{\mathcal L}^k, {\Phi^!_E (C_Z)}}={\H i, {},\Phi_E
({\mathcal L}^k ), {C_Z)}}={\H i, {}, F({\mathcal L}^k ), {C_Z)}}=0
$$
for all $k$ and $i.$

Thus, we
obtain  $\Phi^!_E (C_Z)=0.$ This implies that ${\h \Phi_E (Z), {C_Z}}=0.$
Finally, we get $F(Z)=C_Z [-1]\oplus \Phi_E (Z).$ But we know that
${\h F(Z)[1], {C_Z}}=0.$ Thus, $C_Z=0$ and $f$ is an isomorphism.

2) {\sc Uniqueness}. Suppose there exist two objects $E$ and $E_1$ of
$D^{b}(M\times X)$
such that $\Phi_{E_1}\cong F\cong\Phi_{E_2}.$ Let us consider the complex
(\ref{nob}) over $D^{b}(M\times X)$(see the proof Lemma \ref{obj}).
$$
{\mathcal L}^{ -n}\bt F(R_n)\lto\cdots\lto {\mathcal L}^{ -1}\bt
F(R_1)\lto {\mathcal O}_M \bt F(R_0)
$$

By Lemma \ref{pd1}, there exists a convolution of this complex and
all convolutions  are isomorphic. Let $G\in \db{M\times X}$
be a convolution of the complex (\ref{nob}).
Now consider the following complexes
$$
{\mathcal L}^{ -n}\bt F(R_n)\lto\cdots\lto {\mathcal L}^{ -1}\bt
F(R_1)\lto {\mathcal O}_M \bt F(R_0)\lto E_k
$$
Again by Lemma \ref{pd1}, there exists a unique up to isomorphism convolutions
 of these complexes.

Hence we have the canonical morphisms $G\lto E_1$ and $G\lto E_2.$
Moreover, it has been proved above (see the proof of Lemma \ref{obj})
that $ C_1\oplus E_1\cong G\cong C_2
\oplus E_2$ for large $n,$ where $E_k, C_k$ are
objects of
$\db{M\times X}$ such that ${H}^i (E_k)=0$ for $i\not\in [-a, 0]$ and
${H}^i (C_k)=0$ for $i\not\in [-n-a, -n]$ ($a$ was defined in \ref{boun}).
Thus $E_1$ and $E_2$ are isomorphic.

This completes the proof of Theorem \ref{main} $\Box$

\th{Theorem}\label{eqo} Let $M$ and $X$ be smooth projective varieties.
Suppose $F : \db{M}\lto\db{X}$ is an equivalence. Then there exists a unique
up to isomorphism object
$E\in \db{M\times X}$ such that the functors $F$ and $\Phi_E$ are isomorphic.
\eth
It follows immediately from Theorem \ref{main}

\section{Derived categories of K3 surfaces}
\raz In this chapter we are trying to answer the following question: When
are derived categories of coherent sheaves on two different K3 surfaces over
field $\Bbb{C}$
equivalent?

This question is interesting, because there exists a procedure for recovering
 a variety from its derived category of coherent sheaves if the canonical (or anticanonical)
sheaf is ample. Besides, if $\db{X}\simeq\db{Y}$ and $X$ is a smooth projective
K3 surface, then $Y$ is also a smooth projective K3 surface. This is true,
because the dimension of a variety and Serre functor are invariants of a derived
category.

The following theorem is proved in \cite{BOr}.

\th{Theorem}(see \cite{BOr})\label{rec}
Let $X$ be smooth irreducible projective variety with either ample
canonical or ample anticanonical sheaf. If $D=\db{X}$ is equivalent to $\db{X'}$
for some other smooth algebraic variety, then $X$ is isomorphic to $X'.$
\eth

However, there exist examples of varieties that have equivalent derived
categories, if the canonical sheaf is not ample.
Here we give an explicit description for K3 surfaces with equivalent derived
categories.
\th{Theorem}\label{K3}
Let $S_1$ and $S_2$ be smooth projective K3 surfaces over field $\Bbb{C}.$
Then the derived
categories $\db{S_1}$ and $\db{S_2}$ are equivalent as triangulated categories
iff there exists a Hodge isometry $f_{\tau} : T_{S_1}\stackrel{\sim}{\lto}T_{S_2}$
between the lattices of transcendental cycles of $S_1$ and $S_2.$
\eth
Recall that the lattice of transcendental cycles $T_S$ is the orthogonal
complement to Neron-Severi lattice $N_S$ in $H^2 ( S, {\Bbb{Z}} ).$
{\sf Hodge} isometry means that the one-dimensional subspace $H^{2,0} (S_1)\subset
T_{S_1}\ot\Bbb{C}$ goes to   $H^{2,0} (S_2)\subset T_{S_2}\ot\Bbb{C}.$

Now we need some basic facts about K3 surfaces (see \cite{Mu}).
If $S$ is a K3 surface, then the Todd class $td_S$ of $S$ is equal to $1+2w,$
where $1\in H^0 ( S, {\Bbb{Z}} )$ is the unit element of the cohomology ring
$H^* ( S, {\Bbb{Z}} )$ and $w\in H^4 ( S, {\Bbb{Z}} )$ is the fundamental cocycle
of $S.$ The positive square root $\sqrt{td_S}= 1+w$ lies in $H^* ( S, {\Bbb{Z}} )$
too.

Let $E$ be an object of $\db{S}$ then the Chern character
$$
ch(E) = r(E) + c_1 (E) + \frac{1}{2} (c^2_1 - 2c_2)
$$
belongs to integral cohomology $H^* ( S, {\Bbb{Z}} ).$

For an object $E,$ we put $v(E) = ch(E)\sqrt{td_S}\in H^* ( S, {\Bbb{Z}} )$
and call it the vector associated to $E$ (or Mukai vector).

We can define a symmetric integral bilinear form $(,)$ on $H^* ( S, {\Bbb{Z}} )$
by the rule
$$
(u, u') = rs' + sr' - \alpha\alpha' \in H^4 ( S, {\Bbb{Z}} )\cong{\Bbb{Z}}
$$
for every pair $u=(r, \alpha, s), u'=(r', \alpha', s')\in H^0 ( S, {{\Bbb{Z}}} ) \oplus H^2 ( S, {{\Bbb{Z}}} )\oplus H^4 ( S, {{\Bbb{Z}}} ).$ By $\widetilde{H} ( S, {{\Bbb{Z}}} )$
denote $ H^* ( S, {\Bbb{Z}} )$ with this inner product $(,)$ and call it Mukai
lattice.

For any objects $E$ and $F,$ inner product $(v(E), v(F))$ is equal to the
$H^4$ component  of $ch(E)^{\vee} \cdot ch(F)\cdot td_S.$ Hence, by Riemann-Roch-
Grothendieck theorem, we have
$$
(v(E), v(F)) = \chi (E, F):=\sum_i (-1)^i dim{\E i, {}, E, F}
$$

Let us suppose that $\db{S_1}$ and $\db{S_2}$ are equivalent. By Theorem
\ref{main} there exists an object $E\in\db{S_1\times S_2}$ such that
the functor $\Phi_E$ gives this equivalence.

Now consider the algebraic cycle $Z:=p^*\sqrt{td_{S_1}}\cdot ch(E)\cdot
\pi^*\sqrt{td_{S_2}}$ on the product $S_1\times S_2,$ where $p$ and $\pi$ are
the projections
$$
\begin{array}{ccc}
S_1\times S_2&\stackrel{\pi}{\lto}&S_2\\
\llap{\ss{p}}\da&&\\
S_1&&
\end{array}
$$

It follows from the following lemma that
the cycle $Z$ belongs to integral cohomology $H^* (S_1\times S_2, {\Bbb{Z}}).$
\th{Lemma}\cite{Mu}
For any object $E\in\db{S_1\times S_2}$ the Chern character $ch(E)$ is integral,
which means that it belongs to $H^* (S_1\times S_2, {\Bbb{Z}})$
\eth

The cycle $Z$ defines a homomorphism from integral cohomology of $S_1$ to
integral cohomology of $S_2$:
$$
\begin{array}{cccc}
f:&H^* (S_1, {\Bbb{Z}})&\lto&H^* (S_2, {\Bbb{Z}})\\
&\cup&&\cup\\
&\alpha&\mapsto&\pi_* (Z\cdot p^* (\alpha))
\end{array}
$$
The following proposition is similar to Theorem 4.9 from \cite{Mu}.
\th{Proposition}\label{mu}
If $\Phi_E$ is full and faithful functor from $\db{S_1}$ to $\db{S_2}$
then:

1) $f$ is an isometry between $\widetilde{H} (S_1, {\Bbb{Z}})$  and
$\widetilde{H} (S_2, {\Bbb{Z}}),$

2) the inverse of $f$ is equal to the homomorphism
$$
\begin{array}{cccc}
f':&H^* (S_2, {\Bbb{Z}})&\lto&H^* (S_1, {\Bbb{Z}})\\
&\cup&&\cup\\
&\beta&\mapsto&p_* (Z^{\vee}\cdot \pi^* (\beta))
\end{array}
$$
defined by $Z^{\vee}= p^* \sqrt{td_{S_1}}\cdot ch(E^{\vee})\cdot
\pi^* \sqrt{td_{S_2}},$ where $E^{\vee}:={\pmb R^{\bdot}}{\mathcal
H}om( E, {\mathcal O}_{S_1\times S_2}).$ \eth \pr The left and right
adjoint functors to $\Phi_E$ are:
$$
\Phi_E^* =\Phi_E^! = p_* (E^{\vee}\ot \pi^*(\bdot))[2]
$$
Since $\Phi_E$ is full and faithful, the composition $\Phi_E^*\circ \Phi_E$ is
isomorphic to $id_{\db{S_1}}.$

Functor $id_{\db{S_1}}$ is given by the structure sheaf ${\mathcal
O}_{\Delta}$ of the diagonal $\Delta\subset S_1\times S_1.$

Using the projection formula and Grothendieck-Riemann-Roch theorem,
it can easily be shown that the composition $f'\circ f$ is given by
the cycle $p_1^* \sqrt{td_{S_1}}\cdot ch({\mathcal O}_{\Delta})\cdot
p_2^* \sqrt{td_{S_1}},$ where $p_1, p_2$ are the projections of
$S_1\times S_1$ to the summands. But this cycle is equal to
$\Delta.$

Therefore, $f'\circ f$ is the identity,
and, hence, $f$ is an isomorphism of the lattices, because these lattices
are free abelian groups of the same rank.

Let $\nu_S : S\lto Spec\Bbb{C}$ be the structure morphism of $S.$ Then our
inner product $(\alpha, \alpha')$ on $\widetilde{H} (S, {\Bbb{Z}})$ is equal to
$\nu_*(\alpha^{\vee}  \cdot\alpha').$ Hence, by the projection formula, we have
$$
\begin{array}{rcl}
(\alpha, f(\beta))&=&\nu_{S_2, *}(\alpha^{\vee}\cdot\pi_* (\pi^* \sqrt{td_{S_2}}
\cdot ch(E)\cdot p^* \sqrt{td_{S_1}}\cdot p^* (\beta)))=\\
&=&\nu_{S_2, *}\pi_* (\pi^*(\alpha^{\vee})\cdot p^* (\beta)\cdot ch(E)\cdot
\sqrt{td_{S_1\times S_2}})=\\
&=&\nu_{S_1\times S_2, *} (\pi^* (\alpha^{\vee})\cdot p^* (\beta)\cdot ch(E)
\cdot \sqrt{td_{S_1\times S_2}})
\end{array}
$$
for every $\alpha\in H^* ( S_2, {\Bbb{Z}} ), \beta\in H^* ( S_1, {\Bbb{Z}} ).$
In a similar way, we have
$$
(\beta, f'(\alpha)) =\nu_{S_1\times S_2, *} ( p^* (\beta^{\vee})\cdot
\pi^* (\alpha)\cdot ch(E)^{\vee} \cdot \sqrt{td_{S_1\times S_2}})
$$
Therefore, $(\alpha, f(\beta))=(f'(\alpha), \beta).$ Since $f'\circ f$ is
the identity, we obtain
$$
(f(\alpha), f(\alpha'))=(f'f(\alpha), \alpha')=(\alpha, \alpha')
$$
Thus, $f$ is an isometry. $\Box$

\raz
Consider the isometry $f.$ Since the cycle $Z$ is algebraic, we obtain
two isometries $f_{alg} : -N_{S_1}\bot U\stackrel{\sim}{\lto} -N_{S_2}\bot U$
 and $f_{\tau} : T_{S_1}\stackrel{\sim}{\lto} T_{S_2},$ where $N_{S_1}, N_{S_2}$
are Neron-Severi lattices, and $T_{S_1}, T_{S_2}$ are the lattices of transcendental
cycles. It is clear  $f_{\tau}$ is a Hodge isometry.

Thus we have proved that if the derived  categories of two K3 surfaces are
equivalent, then there exists a Hodge isometry between the lattices of
transcendental cycles.

\raz Let us begin to prove the converse. Suppose we have a Hodge isometry
$$
f_{\tau} : T_{S_2}\stackrel{\sim}{\lto}T_{S_1}
$$
It implies from the following
proposition that we can extend this isometry to Mukai lattices.

\th{Proposition}\cite{Ni}\label{Ni}
Let $\phi_1 , \phi_2 : T\lto H$ be two primitive embedding of a lattice $T$
in an even unimodular lattice $H.$ Assume that the orthogonal complement
$N:=\phi_1 (T)^{\perp}$ in $H$ contains the hyperbolic lattice
$U=\left(\begin{array}{cc}0&1\\1&0\\\end{array}\right)$
as a sublattice.

Then $\phi_1$ and $\phi_2$ are equivalent, that means there exists an
isometry $\gamma$ of $H$ such that $\phi_1 = \gamma\phi_2.$
\eth
We know that the orthogonal complement of $T_S$ in Mukai lattice
$\widetilde{H}( S, {\Bbb{Z}} )$ is isomorphic to $N_S \perp U.$ By Proposition
\ref{Ni}, there exists an isometry
$$
f :  \widetilde{H}( S_2 , {\Bbb{Z}} )\stackrel{\sim}{\lto}
\widetilde{H}( S_1 , {\Bbb{Z}} )
$$
such that $f\bigl|_{T_{S_2}} = f_{\tau}.$

Put $v=f(0,0,1)=(r,l,s)$ and $u=f(1,0,0)=(p,k,q).$

We may assume that $r>1.$ One may do this, because there are two
types of isometries on Mukai lattice that are identity on the lattice of
transcendental cycles. First type is multiplication by Chern character $e^m$
of a line bundle:
$$
\phi_m (r,l,s)=(r,l+rm,s+(m,l)+\frac{r}{2}m^2 )
$$
Second type is the change $r$ and $s$(see \cite{Mu}).
So we can change $f$ in such a way that $r>1$ and $f\bigl|_{T_{S_2}}=f_{\tau}.$

First, note that vector $v\in U\perp -N_{S_1}$ is isotropic, i.e
$(v,v)=0.$ It was proved by Mukai in his brilliant paper \cite{Mu}
that there exists a polarization $A$ on $S_1$ such that the moduli
space ${\mathcal M}_A (v)$ of stable bundles with respect to $A,$
for which vector Mukai is equal to $v,$ is projective smooth K3
surface. Moreover, this moduli space is fine, because there exists
the vector $u\in U\perp -N_{S_1}$ such that $(v,u)=1.$ Therefore we
have a universal vector bundle ${\mathcal E}$ on the product $S_1
\times {\mathcal M}_A (v).$

The universal bundle ${\mathcal E}$ gives the functor
$\Phi_{\mathcal E} : \db{{\mathcal M}_A (v)}\lto\db{S_1}.$

Let us assume that $\Phi_{\mathcal E}$ is an equivalence of  derived
categories. In this case, the cycle
$Z=\pi^*_{S_1}\sqrt{td_{S_1}}\cdot ch({\mathcal E})\cdot
p^*\sqrt{td_{\mathcal M}}$ induces  the Hodge isometry
$$
g : \widetilde{H} ( {\mathcal M}_A (v), {\Bbb{Z}} )
\lto\widetilde{H} ( S_1, {\Bbb{Z}} ),
$$
such that $g(0,0,1)=v=(r,l,s).$
Hence, $f^{-1}\circ g$ is an isometry too, and it sends $(0,0,1)$ to $(0,0,1).$
Therefore $f^{-1}\cdot g$ gives the Hodge isometry between the second cohomologies,
because for a K3 surface $S$
$$
\begin{array}{c}
H^2 ( S, {\Bbb{Z}} )=(0,0,1)^{\perp}\Big/{\Bbb{Z}}(0,0,1).
\end{array}
$$
Consequently, by the strong Torelli theorem (see \cite{Lo}), the
surfaces $S_2$ and ${\mathcal M}_A (v)$ are isomorphic. Hence the
derived categories of $S_1$ and $S_2$ are equivalent.

\raz Thus, to conclude the proof of Theorem \ref{K3}, it remains to
show that the functor $\Phi_{\mathcal E}$ is an equivalence.

First, we show that the functor $\Phi_{\mathcal E}$ is full and
faithful. This is a special case of the following more general
statement, proved in \cite{BO}.

\th{Theorem}\cite{BO}\label{mai}
Let $M$ and $X$ be smooth algebraic varieties and\hfill\\
 $E\in{\db {M\times X}}.$ Then $\Phi_{E}$ is fully faithful functor,
iff the following orthogonality conditions are verified:
$$
\begin{array}{lll}
i) & {\H i, X, \Phi_E({\mathcal O}_{t_1}), {\Phi_{E}({\mathcal
O}_{t_2})}} = 0
& \qquad \mbox{for every }\: i\;\mbox{ and } t_1\ne t_2.\\
&&\\
ii) & {\H 0, X, \Phi_E({\mathcal O}_t), {\Phi_E({\mathcal O}_t)}} = k,&\\
&&\\
& {\H i, X, \Phi_E({\mathcal O}_t), {\Phi_E({\mathcal O}_t)}} = 0 ,
& \qquad \mbox{ for }i\notin [0, dim M].
\end{array}
$$
Here $t,$ $t_{1},$ $t_{2}$ are points of $M,$ ${\mathcal O}_{t_{i}}$
are corresponding skyscraper sheaves. \eth

In our case, $\Phi_{\mathcal E} ({\mathcal O}_t)=E_t,$ where $E_t$
is stable sheaf with respect to the polarization $A$ on $S_1$ for
which $v(E_t)=v.$ All these sheaves are simple and ${\E i, {}, E_t ,
{E_t}}=0$ for $i\not\in [0,2].$ This implies that condition 2) of
Theorem \ref{mai} is fulfilled.

All $E_t$ are stable sheaves, hence ${\h E_{t_1}, {E_{t_2}}}=0.$ Further,
by Serre duality ${\E 2, {}, E_{t_1}, {E_{t_2}}}=0.$ Finally, since the vector
$v$ is isotropic, we obtain ${\E 1, {}, E_{t_1}, {E_{t_2}}}=0.$

This yields that $\Phi_{\mathcal E}$ is full and faithful. As  our
situation is not symmetric (a priori), it is not clear whether  the
adjoint functor to $\Phi_{\mathcal E}$ is also full and faithful.
Some additional reasoning is needed.

\th{Theorem}\label{equi} In the above notations, the functor
$\Phi_{\mathcal E}: \db{{\mathcal M}_A (v)}\lto\db{S_1}$
 is an equivalence.
\eth \pr Assume the converse, i.e. $\Phi_{\mathcal E}$ is not an
equivalence, then,
 since the functor $\Phi_{\mathcal E}$ is full and faithful, there
exists an object $C\in \db{S_1}$ such that $\Phi^*_{\mathcal E}
(C)=0.$ By Proposition \ref{mu}, the functor $\Phi_{\mathcal E}$
induces the isometry $f$ on the Mukai lattices, hence the Mukai
vector $v(C)$ is equal to $0.$

Object $C$ satisfies the conditions ${\H i, {}, C, {E_t}}=0$ for
every $i$ and all $t\in {\mathcal M}_A (v),$ where $E_t$ are stable
bundles on $S_1$ with the Mukai vector $v.$

Denote by $H^i (C)$ the cohomology sheaves of the object $C.$
There is a spectral sequence which converges to ${\H i, {}, C, {E_t}}$
\begin{equation}\label{seq}
E^{p,q}_2 = {\E p, {}, H^{-q} (C), {E_t}} \Longrightarrow {\H {p+q}, {}, C, {E_t}}
\end{equation}
It is depicted in the following diagram

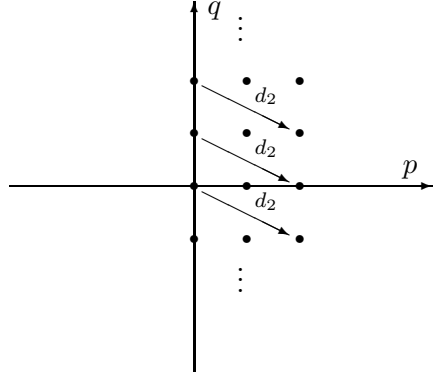
\begin{picture}(400,160)
\put(200,10){\vector(0,1){140}}
\put(130,80){\vector(1,0){160}}
\put(201,145){$q$}
\put(275,85){$p$}

\put(200,120){\circle*{3}}
\put(200,100){\circle*{3}}
\put(200,60){\circle*{3}}
\put(200,80){\circle*{3}}
\put(240,120){\circle*{3}}
\put(240,100){\circle*{3}}
\put(240,80){\circle*{3}}
\put(240,60){\circle*{3}}
\put(220,120){\circle*{3}}
\put(220,100){\circle*{3}}
\put(220,80){\circle*{3}}
\put(220,60){\circle*{3}}
\put(212,135){$\vdots$}
\put(212,40){$\vdots$}
\put(203,78){\vector(2,-1){33}}
\put(203,98){\vector(2,-1){33}}
\put(203,118){\vector(2,-1){33}}
\put(219,72){\scriptsize{$d_2$}}
\put(219,92){\scriptsize{$d_2$}}
\put(219,112){\scriptsize{$d_2$}}

\end{picture}

We can see that ${\E 1, {}, H^q (C), {E_t}}=0$ for every $q$ and all $t,$
and every morphism $d_2$ is an isomorphism.

To prove the theorem, we need the following lemma.

\th{Lemma}\label{GN}
Let $G$ be a sheaf on K3 surface $S_1$ such that ${\E 1, {}, G,
{E_t}}=0$ for all $t.$ Then there exists an exact sequence
$$
0\lto G_1 \lto G\lto G_2 \lto 0
$$
that satisfies the following conditions:
$$
\begin{array}{llll}
1)\; {\E i, {}, G_1, {E_t}}=0 \quad \mbox{ for every }\; i\ne 2,\;\mbox{and}
& {\E 2, {}, G_1, {E_t}}\cong{\E 2, {}, G, {E_t}}       \\
2)\; {\E i, {}, G_2, {E_t}}=0 \quad \mbox{ for every } \; i\ne 0,\;\mbox{and}
& {\H {}, {}, G_2, {E_t}}\cong{\H {}, {}, G, {E_t}}
\end{array}
$$
and $p_A (G_2) < p_A (G) < p_A (G_1)$ (where $p_A(F)$ is a Gieseker slope, i.e., a polynomial such that $p_A(F)(n)=\chi(F(nA))/r(F).$)
\eth
\pr
Firstly, there is a short exact sequence
$$
0\lto T\lto G\lto \widetilde{G}\lto 0,
$$
where $T$ is a torsion sheaf, and $\widetilde{G}$ is torsion free.

Secondly, there is a Harder-Narasimhan filtration $0=I_0\subset ...\subset I_n
=\widetilde{G}$
for $\widetilde{G}$
such that the successive quotients $I_i / I_{i-1}$ are $A$-semistable, and
$p_A (I_i / I_{i-1})>p_A( I_j / I_{j-1} )$ for $i<j.$

Now, combining $T$ and the members of the filtration for which $p_A (I_i
/ I_{i-1})>p_A (E_t)$ (resp. $=,$ $<$) to  one, we obtain the 3-member
 filtration on $G$
$$
0=J_0 \subset J_1 \subset J_2 \subset J_3 =G.
$$
Let $K_i$ be the quotients sheaves $J_i / J_{i-1}.$ We have
$$
p_A (K_1) > p_A (K_2)=p_A (E_t) > p_A (K_3)
$$
(we suppose, if  needed, $p_A (T)= +\infty$).

Moreover, it follows from  stability  of $E_t$ that
$$
{\h K_1, {E_t}}=0 \qquad \mbox{ and }\qquad {\E 2, {}, K_3, {E_t}}=0
$$
Combining this with the assumption that ${\E 1, {}, G, {E_t}}=0,$ we get
${\E 1, {}, K_2, {E_t}}=0.$

To prove the lemma it remains to show that $K_2 =0.$

Note that $K_2$ is $A$-semistable. Hence there is a Jordan-H\"older filtration
for $K_2$ such that  the successive quotients are $A$-stable. The number of
the quotients is finite. Therefore we can take $t_0$ such that
$$
{\h K_2, {E_{t_0}}}=0 \qquad \mbox{ and }\qquad {\E 2, {}, K_2, {E_{t_0}}}=0
$$
Consequently, $\chi ( v(K_2), v(E_t) )=0.$ Thus, as
${\E 1, {}, K_2, {E_t}}=0$ for all $t,$ we obtain
${\E i, {}, K_2, {E_t}}=0$ for every $i$ and all $t.$

Further, let us consider $\Phi^*_{\mathcal E} (K_2).$ We have
$$
{\H \bdot , {},  \Phi^*_{\mathcal E}(K_2), {{\mathcal O}_t}}\cong
{\H \bdot , {}, K_2, {E_t}}=0,
$$
This implies $\Phi^*_{\mathcal E} (K_2)=0.$ Hence $v(K_2)=0,$
because $f$ is an isometry. And, finally, $K_2=0.$ The lemma is
proved. $\Box$

Let us return to the theorem. The object $C$ possesses at least two
 non-zero consequent cohomology
sheaves $H^p (C)$ and $H^{p+1} (C)$ . They satisfy the
condition of Lemma \ref{GN}  Hence there exist decompositions with conditions 1),2):
$$
0 \lto H^p_1 \lto H^p (C) \lto H^p_2 \lto 0 \quad\mbox{and}\quad
0 \lto H^{p+1}_1 \lto H^{p+1} (C) \lto H^{p+1}_2 \lto 0
$$
Now consider the canonical morphism $H^{p+1} (C) \lto H^p (C)[2].$ It
induces the morphism $ s : H^{p+1}_1 \lto H^p_2 [2].$ By $Z$ denote a cone of
$s.$

 Since $d_2$ of the spectral sequence (\ref{seq}) is
an isomorphism, we obtain
$$
{\H \bdot, {}, Z, {E_t}}=0 \qquad \mbox{ for all } t.
$$
Consequently, we have $\Phi^*_{\mathcal E} (Z)=0.$ On the other
hand, we know that $p_A (H^{p+1}_1)>p_A (E_t)>p_A (H^p_2).$
 Therefore $v(Z)\ne 0.$
This contradiction proves the theorem. $\Box$

There exists the another version of Theorem \ref{K3}
\th{Theorem}
Let $S_1$ and $S_2$ be smooth projective K3 surfaces over field $\Bbb{C}.$
Then the derived
categories $\db{S_1}$ and $\db{S_2}$ are equivalent as triangulated categories
iff there exists a Hodge isometry $f: \widetilde{H} ( S_1, {\Bbb{Z}} )\stackrel{\sim}\lto
\widetilde{H} ( S_2, {\Bbb{Z}} )$
between the Mukai lattices of $S_1$ and $S_2.$
\eth
Here the `{\sf Hodge} isometry' means that the one-dimensional subspace $H^{2,0} (S_1)\subset
\widetilde{H} ( S_1, {\Bbb{Z}} )\ot\Bbb{C}$ goes to   $H^{2,0} (S_2)\subset
\widetilde{H} ( S_2, {\Bbb{Z}} )\ot\Bbb{C}.$

\sec{Appendix.}

 The facts, collected in this appendix, are not new; they are known. However,
not having a good reference, we regard it necessary to give a proof for the statement, which is used in the main text. We exploit the technique from \cite{IM}.

Let $X$ be a smooth projective variety and $L$ be a very ample invertible sheaf on $X$ such that ${\rm H}^i ( X , L^k ) =0$ for any $k>0$ , when $i\ne 0.$
Denote by $A$ the coordinate algebra for $X$ with respect to $L,$ i.e. $A = \bigoplus\limits^{\infty}_{k=0} {\rm H}^0 ( X , L^k ).$

Now consider the variety $X^n.$ First, we introduce some notations. Define subvarieties $\Delta^{(n)}_{(i_1 , ..., i_k )(i_{k+1}, ..., i_m )} \subset X^n$ by the following rule:
$$
\Delta^{(n)}_{(i_1 , ..., i_k )(i_{k+1}, ..., i_m )} :=\{ (x_1 , ..., x_n ) | x_{i_1}=\cdots = x_{i_k} ; x_{i_{k+1}}=\cdots = x_{m} \}
$$
By $S^{(n)}_i$ denote $\Delta^{(n)}_{(n,..., i)}.$ It is clear that $S^{(n)}_i \cong X^i.$

Further, let $T^{(n)}_i := \bigcup\limits^{i-1}_{k=1}
\Delta^{(n)}_{(n,..., i)(k, k-1)}$ (note that $T^{(n)}_1$ and
$T^{(n)}_2$ are empty) and let $\Sigma^{(n)} :=
\bigcup\limits^{n}_{k=1} \Delta^{(n)}_{(k, k-1)}.$ We see that
$T^{(n)}_i \subset S^{(n)}_i.$ Denote by ${\mathcal I}^{(n)}_i$ the
kernel of the restriction map ${\mathcal O}_{S^{(n)}_i}\lto
{\mathcal O}_{T^{(n)}_i}\lto 0.$

Using induction by $n,$ it can easily be checked that the following complex on $X^n$
$$
P^{\bdot}_n  : 0\lto J_{\Sigma^{(n)}} \lto {\mathcal I}^{(n)}_n \lto
{\mathcal I}^{(n)}_{n-1} \lto \cdots \lto {\mathcal I}^{(n)}_2 \lto
{\mathcal I}^{(n)}_1 \lto 0
$$
is exact. (Note that ${\mathcal I}^{(n)}_1 = {\mathcal
O}_{\Delta^{(n)}_{n,...,1}}$ and ${\mathcal I}^{(n)}_2 = {\mathcal
O}_{\Delta^{(n)}_{n,...,2}}$). For example, for $n=2$ this complex
is a short exact sequence on $X\times X$:
$$
P^{\bdot}_2 : 0\lto J_{\Delta}\lto {\mathcal O}_{X\times X}\lto
{\mathcal O}_{\Delta}\lto 0
$$

Denote by $\pi^{(n)}_i$ the projection of $X^n$ onto $i^{th}$ component, and by $\pi^{(n)}_{ij}$ denote the projection of $X^n$ onto the product of $i^{th}$ and $j^{th}$ components.

Let $B_n := {\rm H}^0 ( X^n , J_{\Sigma^{(n)}}\ot (L\bt \cdots \bt
L))$ and let $R_{n-1} := R^{0} \pi^{(n)}_{1 *} (J_{\Sigma^{(n)}}\ot
({\mathcal O}\bt L\bt \cdots \bt L)).$ \ap{Proposition} Let $L$ be a
very ample invertible sheaf on $X$ as above. Suppose that for any
$m$ such that $1< m \le n+dim X+2$ the following conditions hold:
$$
\begin{array}{lll}
a)& {\rm H}^i ( X^m , J_{\Sigma^{(m)}}\ot (L\bt \cdots \bt L))=0&\; \mbox{for} \quad i\ne 0\\
b)& R^{i} \pi^{(m)}_{1 *} (J_{\Sigma^{(m)}}\ot ({\mathcal O}\bt L\bt \cdots \bt L))=0&\; \mbox{for}\quad i\ne 0\\
c)& R^{i} \pi^{(m)}_{1m *} (J_{\Sigma^{(m)}}\ot ({\mathcal O}\bt
L\bt \cdots \bt L\bt {\mathcal O}))=0&\; \mbox{for} \quad i\ne 0
\end{array}
$$
Then we have:

1) algebra $A$ is n-Koszul, i.e the sequence
$$
B_n \ot_k A\lto B_{n-1}\ot_k A\lto \cdots\lto B_1\ot_k A\lto A\lto k\lto 0
$$
is exact;

2)   the following
complexes on $X$:
$$
A_{k-n}\ot R_n \lto A_{k-n+1}\ot R_{n-1}\lto\cdots\lto A_{k-1}\ot
R_1 \lto A_k \ot R_0 \lto {\mathcal L}^{k}\lto 0
$$
are exact for any $k\ge 0$
(if $ k-i<0,$ then $A_{k-i}=0$ by definition);

3)  the  complex
$$
{ L}^{ -n}\bt R_n\lto\cdots\lto {L}^{ -1}\bt R_1\lto {\mathcal O}_M
\bt R_0 \lto {\mathcal O}_{\Delta}
$$
gives n-resolution of the diagonal on $X\times X,$ i.e. it is exact.
\eap
\pr
1) First, note that
$$
{\rm H}^i ( X^m , {\mathcal I}^{(m)}_k \ot (L\bt \cdots \bt L))={\rm
H}^i ( X^{k-1} , J_{\Sigma^{(k-1)}}\ot (L\bt \cdots \bt L))\ot
A_{m-k+1}
$$
By condition a), they are trivial for $i\ne 0.$

Consider the complexes $P^{\bdot}_m \ot (L\bt\cdots \bt L)$ for $m\le n+dim X +1.$ Applying the functor ${\rm H}^0$ to these complexes and using condition a), we get the exact sequences:
$$
0\lto B_m \lto B_{m-1}\ot_k A_1 \lto \cdots\lto B_1\ot_k A_{m-1} \lto A_m \lto 0
$$
for $m\le n+dimX+1.$

Now put $m=n+dimX+1.$ Denote by $W^{\bdot}_m$ the complex
$$
{\mathcal I}^{(m)}_m \lto {\mathcal I}^{(m)}_{m-1} \lto \cdots \lto
{\mathcal I}^{(m)}_2 \lto {\mathcal I}^{(m)}_1 \lto 0
$$
Take the complex $W^{\bdot}_m \ot (L\bt\cdots\bt L\bt L^i )$ and apply functor
${\rm H}^0$ to it. We obtain the following sequence:
$$
 B_{m-1}\ot_k A_i \lto B_{m-2}\ot_k A_{i+1} \lto\cdots\lto B_1\ot_k A_{m-1} \lto A_m \lto 0
$$
Its cohomologies are ${\rm H}^j ( X^m , J_{\Sigma^{(m)}}\ot (L\bt \cdots \bt L\bt L^i)).$ It follows from condition b) that
$$
{\rm H}^j ( X^m , J_{\Sigma^{(m)}}\ot (L\bt \cdots \bt L\bt
L^i))={\rm H}^j ( X,  R^0 \pi^{(m)}_{m*}(J_{\Sigma^{(m)}}\ot (L\bt
\cdots \bt L\bt {\mathcal O}))\ot L^i ).
$$
Hence they are trivial for $j> dimX.$
Consequently, we have the exact sequences:
$$
B_n \ot_k A_{m-n+i-1}\lto B_{n-1}\ot_k A_{m-n+i} \lto \cdots\lto B_1\ot_k A_{m+i-2}\lto A_{m+i-1}
$$
for $i\ge 1.$ And for $i\le 1$ the exactness was proved above.
Thus, algebra $A$ is n-Koszul.

2) The proof is the same as for 1). We have isomorphisms
$$
R^i \pi^{(m)}_{1*}({\mathcal I}^{(m)}_k \ot ({\mathcal O}\bt L\bt
\cdots \bt L))\cong R^i \pi^{(k-1)}_{1*}( J_{\Sigma^{(k-1)}}\ot
({\mathcal O}\bt L\bt \cdots \bt L))\ot A_{m-k+1}
$$
Applying functor $R^0 \pi^{(m)}_{1*}$ to the complexes
$P^{\bdot}_{m} \ot ({\mathcal O}\bt L\bt \cdots \bt L))$ for $m\le
n+dimX+2,$ we obtain the exact complexes on $X$
$$
0\lto R_{m-1} \lto A_{1}\ot R_{m-2}\lto\cdots\lto A_{m-2}\ot R_1
\lto A_{m-1} \ot R_0 \lto {\mathcal L}^{m-1}\lto 0
$$
for $m\le n+dimX+2.$

Put $m=n+dimX+2.$ Applying functor $R^0 \pi^{(m)}_{1*}$ to the
complex $W^{\bdot}_m \ot ({\mathcal O}\bt L\bt \cdots \bt L \bt
L^i)),$ we get the complex
$$
 A_{i}\ot R_{m-2}\lto\cdots\lto A_{m+i-3}\ot R_1 \lto
A_{m+i-2} \ot R_0 \lto {\mathcal L}^{m+i-2}\lto 0
$$
The cohomologies of this complex are
$$
R^{j} \pi^{(m)}_{1 *} (J_{\Sigma^{(m)}}\ot ({\mathcal O}\bt L\bt
\cdots \bt L \bt L^i))\cong R^{j}p_{1*}(R^{0} \pi^{(m)}_{1m *}
(J_{\Sigma^{(m)}}\ot ({\mathcal O}\bt L\bt \cdots \bt L\bt {\mathcal
O}))\ot ({\mathcal O}\bt L^i))
$$
They are trivial for $j> dimX.$
Thus, the sequences
$$
A_{k-n}\ot R_n \lto A_{k-n+1}\ot R_{n-1}\lto\cdots\lto A_{k-1}\ot
R_1 \lto A_k \ot R_0 \lto {\mathcal L}^{k}\lto 0
$$
are exact for all $k\ge 0.$

3) Consider the complex $W^{\bdot}_{n+2}\ot ({\mathcal O}\bt L\bt
\cdots \bt L \bt L^{-i}).$ Applying the functor $R^0
\pi^{(n+2)}_{1(n+2)*}$ to it, we obtain the following complex on
$X\times X$:
$$
{L}^{-n}\bt R_n\lto\cdots\lto {L}^{-1}\bt R_1\lto {\mathcal O}_M \bt
R_0 \lto {\mathcal O}_{\Delta}
$$
By condition c), it is exact.
This finishes the proof.

Note that for any ample invertible sheaf $L$ we can find $j$ such that for the
sheaf $L^j$ the conditions a),b),c) are fulfilled.

\end{document}